\documentclass[12pt]{article}
\usepackage{graphicx}
\usepackage{amssymb, amsfonts, amsthm, amsmath}

\theoremstyle{plain}
\newtheorem{theo}{Theorem}

\newtheorem{cor}[theo]{Corollary}

\newtheorem{lem}[theo]{Lemma}

\theoremstyle{definition}

\newtheorem{defi}[theo]{Definition}

\newcommand{\defn}[1]{{\bf #1}}

\newcommand{\cA}{{\cal A}}\newcommand{\cB}{{\cal B}}

\newcommand{\cR}{{\cal R}}

\newcommand{\E}{{\mathbb E}}

\newcommand{\N}{{\mathbb N}}\renewcommand{\P}{{\mathbb P}}
 \newcommand{\R}{{\mathbb R}}

\newcommand{\1}{\bold 1}

\newcommand{\V}{{\text{\rm Var}}}
\newcommand{\Cov}{{\text{\rm Cov}}}
\newcommand{\Corr}{{\text{\rm Corr}}}

\begin{document}

\title{Estimation of ordinal pattern probabilities in fractional Brownian motion}

\author{Mathieu Sinn, Karsten Keller\\
Institute of Mathematics, University of L\"ubeck}

\maketitle

\begin{abstract}
For equidistant discretizations of fractional Brownian
motion (fBm), the probabilities of ordinal patterns of order $d=2$
are monotonically related to the Hurst parameter $H$. 
By plugging the sample relative frequency of 
those patterns indicating changes between up and down
into the monotonic relation to $H$,
one obtains the Zero Crossing (ZC) estimator of the Hurst parameter
which has found considerable attention in mathematical and applied
research.

In this paper, we generally discuss
the estimation of ordinal pattern probabilities in fBm.
As it turns out, according to the sufficiency principle, 
for ordinal patterns of order $d=2$ any reasonable estimator is an affine
functional of the sample relative frequency of changes.
We establish strong consistency of the estimators 
and show
them to be asymptotically normal for $H<\frac{3}{4}$. 
Further,
we derive confidence intervals for the Hurst parameter.
Simulation studies show that the ZC estimator has larger variance
but less bias than the HEAF estimator of the Hurst parameter.
\end{abstract}


\section{Introduction}\label{intro}
Probabilities of ordinal patterns in equidistant dis\-cretizations
of fractional Brownian motion (fBm) have first been analyzed by
Bandt and Shiha \cite{BandtShiha05}. Ordinal patterns represent
the rank order of successive equidistant values of a discrete time
series. For example, for three successive values there are six
different possible outcomes of the rank order, which we call 
ordinal patterns of order $d=2$. 

As Bandt and Shiha have shown, for equidistant discretizations of fBm 
the distribution of ordinal patterns
is stationary and does not
depend on the particular sampling interval length.
Further,
the probabilities of ordinal
patterns of order $d=2$ are all monotonically related to the Hurst
parameter, with the probability of
ordinal patterns indicating changes from up to down and from
down to up, respectively, being strictly monotonically decreasing in $H$.

The estimator of the Hurst parameter
obtained by plugging the sample relative frequency of changes
between up and down
into the monotonic functional
relation to $H$ 
has been known for some time, running under the label
`Zero Crossing' (ZC) estimator since changes
between up and down correspond to
zero crossings of the first order differences. 

Note that more
generally, but not focussing on fBm, Kedem \cite{Kedem94}
has considered the estimation of parts of the
autocorrelation structure of a stationary Gaussian process by
counting (higher order) zero crossings, that is, zero crossings of
the first (or higher) order differences.

The ZC estimator is asymptotically normally distributed for
$H<\frac{3}{4}$, a result essentially due to Ho and Sun \cite{HoSun87}, who proved
that the sample relative frequency of changes is asymptotically normally
distributed in this case.
Coeurjolly \cite{Coeurjolly00} has
resumed properties of the ZC estimator including strong consistency. 
In the sequel, the applicability of the ZC estimator has been examined
by Markovi\'c and Koch \cite{MarkovicKoch05} and Shi et {\it al.} \cite{Shietal05}, 
including simulation studies as well as the
application to hydrological and meteorological data.

The distribution of ordinal patterns is the base of ordinal time
series analysis, being a new fast, robust and flexible approach to
the investigation of large and complex time series (see Bandt \cite{Bandt04}
and Keller et {\it al.} \cite{KellerSinnEmonds07,KellerLaufferSinn05}). From the viewpoint
of ordinal time series analysis, the estimation of ordinal pattern
probabilities in stochastic processes is of special interest. 

This paper is structured as follows: 
Sec.~\ref{estopp} is devoted to a general discussion of the
estimation of ordinal pattern probabilities in fBm where
we consider estimators
based on counting the occurence of ordinal patterns in realizations. 
Obviously, the sample relative frequency of ordinal patterns is
an unbiased estimate of the probability of the patterns.
As Theorem \ref{RaoBlackwell} shows,
averaging the sample relative frequencies of ordinal patterns and of their
`time' and `spatial' reversals, one obtains estimates 
which are strictly more concentrated in 
convex order. Hence, we restrict our following considerations 
to such `reasonable' estimators.
Notice that Theorem \ref{RaoBlackwell} has important consequences for
the estimation of functionals of ordinal pattern distributions such
as the permutation entropy (see \cite{Bandt04}, \cite{KellerLaufferSinn05}).

For ordinal patterns of order $d=2$
one obtains two different estimators which can both
be expressed as affine functionals of the sample relative
frequency of changes.
By Lemma \ref{strongconsist}
we establish strong consistency 
and asymptotical unbiasedness for the estimators
of bounded continuous functionals of ordinal pattern probabilities.
Theorem \ref{asnormtheo} states that
estimators of ordinal pattern probabilities are asymptotically normal
for $H<\frac{3}{4}$.

In Sec.~\ref{var}, we consider the estimation of
the probability of a change by the sample relative frequency of changes
in some detail. We give a 
formula for the precise numerical evaluation
of the variance of the sample relative frequency of changes as well as 
asymptotically equivalent expressions.
Based on these results, confidence intervals for the Hurst parameter 
are provided in Sec.~\ref{conf}.

In Sec.~\ref{sim}, we compare the previous results
to findings for
simulations of fBm. As it turns out, the confidence intervals obtained from
the ZC estimator in Sec.~\ref{conf} cover the true unkown value of the Hurst parameter
for about $95$ per cent of the cases, even for small sample sizes and for
the values of the Hurst parameter larger than $\frac{3}{4}$, where
asymptotic normality of the ZC estimator is not necessarily expected to hold.
Compared to the HEAF estimator which estimates the Hurst parameter by plugging
the sample autocovariance into a monotonic functional
relation to $H$, the ZC estimator has larger variance but much less bias,
in particular for small sample sizes.

Notice that the results given in this work for the
increments of fBm similarly apply to
FARIMA($0$,$d$,$0$) processes. Furthermore, the statistical properties of the estimates 
of ordinal pattern probabilities do not only
apply to fBm, but also to montonic transformations of fBm. In particular,
the estimates of ordinal pattern probabilities are invariant with respect
to (unknown) non-linear monotonic transformations of processes.

\section{Estimating ordinal pattern probabilities }\label{estopp}

\subsection{General aspects}\label{estordp}
Let $(\Omega,\cA)$ be a
measurable space 
equipped with a family
$(\P_\vartheta)_{\vartheta\in\Theta}$ of probability measures with
$\Theta$ non-empty.
The subscript $\vartheta$ indicates that a
corresponding quantity ($\E_{\vartheta}$, $\V_{\vartheta} $, etc.) is taken
with respect to $\P_\vartheta$. For the special case of fractional Brownian
motion where $\Theta={]0,1]}$, we write $H$ instead of $\vartheta$.

Let $(X_k)_{k\in\N_0}$ be a given real-valued stochastic process 
defined on $(\Omega,\cA)$. Define the process of increments $(Y_k)_{k\in\N}$ by
$Y_k:=X_k-X_{k-1}$ for $k\in\N$. For $d\in\N$ let $S_d$ denote
the set of permutations of $\{0,1,\ldots,d\}$.

\begin{defi}
For $d\in\N$ let the mapping $\pi:\R^{d+1}\to S_d$ be
defined by
\begin{eqnarray*}
\pi( (x_0,x_1,\ldots,x_d)  ) &=& \left(
\begin{array}{ccccc}
0 & 1 & 2 & \ldots & d \\
r_0 & r_1 & r_2 & \ldots & r_d
\end{array}
 \right) \ =: \ (r_0,r_1,\ldots,r_d)
\end{eqnarray*}
for $ (x_0,x_1,\ldots,x_d) \in\R^{d+1}$ with the
permutation $(r_0,r_1,\ldots,r_d)$ of $\{0,1,\ldots,d\}$ satisfying 
$x_{d-r_0}\geq
x_{d-r_1}\geq \ldots \geq x_{d-r_d}$,
and $r_{l-1}>r_l$ if $x_{d-r_{l-1}}=x_{d-r_l}$
for $l=1,2,\ldots,d$. For $k\in\N_0$ we call
\begin{eqnarray*}
\Pi_d(k) &:=& \pi( (X_k, X_{k+1}, \ldots, X_{k+d} ))
\end{eqnarray*}
the \textit{(random) ordinal pattern} of \textit{order} $d$
at \textit{time} $k$.
\hfill $\square$ 
\end{defi}

The permutation $\pi( (x_0,x_1,\ldots,x_d)  )$ describes the rank order of the
values $x_0,x_1,\ldots,x_d$, where in case $ x_{d-i} =
x_{d-j} $ for $0\leq i < j \leq d$ the `earlier'
$x_{d-j} $ is ranked higher than $x_{d-i}$.

Notice that if $(X_k)_{k\in\N_0}$ is \textit{pairwise distinct},
that is, $X_i\neq X_j$ $\P_{\vartheta} $-a.s. for all
$i,j\in \N_0$ with $i\neq j$ and $\vartheta\in\Theta$,
then, apart
from sets with probability zero for all $\vartheta\in\Theta$, $\Pi_d(k)$ generates the same $\sigma$-algebra 
as the rank vector $(R_0,R_1,\ldots,R_d)$ of $(X_k, X_{k+1}, \ldots, X_{k+d} )$ given by
$R_j=\sum_{i=0}^d\1_{\{X_{k+j}\geq X_{k+i}\}}$ for $j=0,1,\ldots,d$ 
(see \cite{Lehmann86}, p.~286).

\paragraph{Stationarity.} 
For $\vartheta\in\Theta$ let $\stackrel{ \P_{\vartheta}  }{=}$ denote
equality in distribution with respect to $ \P_{\vartheta} $. 
A stochastic process  $(Z_k)_{k\in T}$ defined on $(\Omega,\cA)$ 
for $T=\N$ or $T=\N_0:=\N\cup\{0\} $ is called
\textit{stationary}
iff with respect to each $\vartheta\in\Theta$
\begin{eqnarray*}
\hspace{10mm} \big( Z_{k_1},Z_{k_2},\ldots,Z_{k_n}\big) 
&\stackrel{ \P_{\vartheta} }{=}& \big( Z_{k_1+l},Z_{k_2+l},\ldots,Z_{k_n+l}\big) 
\end{eqnarray*}
for all $k_1,k_2,\ldots,k_n\in T$ with $n\in\N$ and for all $l\in\N$. 

In fact, $\Pi_d(k)$ only depends on $(Y_{k+1},Y_{k+2},\ldots,Y_{k+d})$ for all $d\in\N$ and $k\in\N_0$.
In particular, let
\begin{eqnarray*}
\widetilde{\pi}( (y_1, y_2, \ldots, y_d)) &=& (r_0,r_1,\ldots,r_d)
\end{eqnarray*}
be the unique permutation of $\{0,1,\ldots,d\}$ for $(y_1, y_2, \ldots, y_d)\in\R^d$ such that
\begin{eqnarray}\label{patincr}
\sum_{j=1}^{d-r_{0}}y_{j} \geq \sum_{j=1}^{d-r_1}y_{j}\geq  \ldots
\geq \sum_{j=1}^{d-r_{d-1}}y_{j} \geq  \sum_{j=1}^{d-r_{d}}y_{j},
\end{eqnarray}
and $r_{l-1}>r_l$ if 
$\sum_{j=1}^{d-r_{l-1}}y_{j} = \sum_{j=1}^{d-r_{l}}y_{j}$
for $l=1,2,\ldots,d$. Obviously,
\begin{eqnarray*}
\pi((x_0,x_1,\ldots,x_d)) &=& \widetilde{\pi}(( x_1-x_0, x_2-x_1, \ldots, x_d-x_{d-1}  ))
\end{eqnarray*}
for all $ (x_0,x_1,\ldots,x_d) \in\R^{d+1}$, and hence
$ \Pi_d(k) = \widetilde{\pi}((Y_{k+1},Y_{k+2},\ldots,Y_{k+d}))$
for all $k\in\N_0$. This immediately yields the following statement.

\begin{cor}\label{statcor}
If $(Y_k)_{k\in\N}$ is stationary 
then $\big( \Pi_d(k)\big)_{k\in\N_0} $ is
stationary.
\end{cor}

\begin{figure}[h]
\begin{picture}(0,80)
\put(20,20){\includegraphics[width=3.5cm]{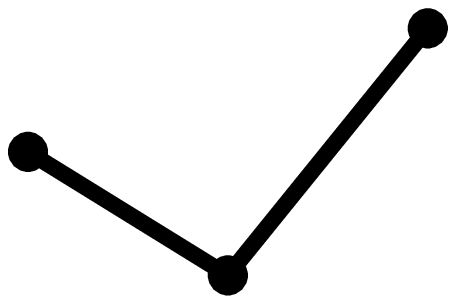}}
\put(120,20){\includegraphics[width=3.5cm]{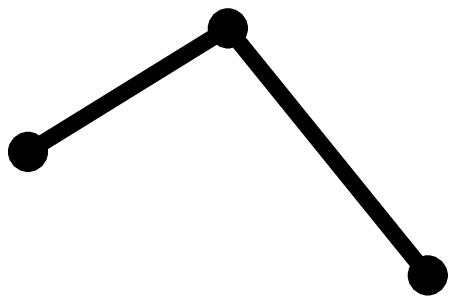}}
\put(220,20){\includegraphics[width=3.5cm]{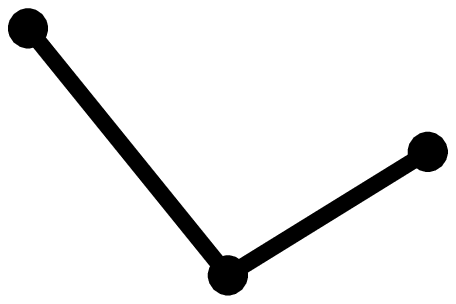}}
\put(320,20){\includegraphics[width=3.5cm]{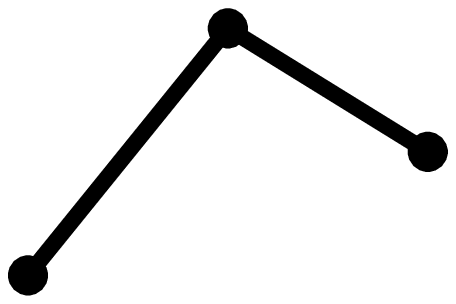}}
\put(40,0){\text{$r=(0,2,1)$}}
\put(130,0){\text{$\alpha(r)=(1,2,0)$}}
\put(230,0){\text{$\beta(r)=(2,0,1)$}}
\put(330,0){\text{$\beta\circ\alpha(r)=(1,0,2)$}}
\end{picture}
\caption{{The ordinal pattern $r=(0,2,1)$, its spatial reversal $\alpha(r)$,
its time reversal $\beta(r)$, and its spatial and time reversal
$\beta\circ\alpha(r)$. }}\label{fig2}
\end{figure}

\paragraph{Space and time symmetry.} For $d\in\N$ let the mappings
$\alpha,\beta$ from $S_d$ onto $S_d$ be defined by
\begin{eqnarray}\label{alphabeta}
\alpha(r) := (r_d,r_{d-1},\ldots,r_0) , \ \ \ \beta(r) :=
(d-r_0,d-r_{1},\ldots,d-r_d)
\end{eqnarray}
for $r=(r_0,r_1,\ldots,r_d)\in S_d$. Geometrically, $\alpha(r)$ and
$\beta(r)$ can be seen as the spatial and time reversal of $r$,
respectively (see Figure \ref{fig2}). Let the set
$\overline{r}$ be defined by
\begin{eqnarray}\label{rclos}
\overline{r} &:=& \{r, \alpha(r), \beta(r),
\beta\circ\alpha(r)\}
\end{eqnarray}
with $\circ$ denoting the usual composition of functions.
As $\alpha\circ\beta(r)=\beta\circ\alpha(r)$ and
$\alpha\circ\alpha(r)=\beta\circ\beta(r)=r$, one has
$\alpha(\overline{r})=\beta(\overline{r})=\overline{r}$.
Clearly, if $s\in\overline{r}$ for $r,s\in S_d$, 
then $\overline{s}=\overline{r}$. This
provides a division of each $S_d$ into classes, consisting of 2 or 4
elements. For $d=1$ the only class is $S_d=\{(0,1),(1,0)\}$, for
$d=2$ there are the two classes $\{(0,1,2),(2,1,0)\}$ and
$\{(0,2,1),(2,0,1),(1,2,0),(1,0,2)\}$, and for $d=3$ there are 8
classes. Note that for $d\geq 3$ classes of both $2$ and $4$
elements are possible.

For $d\in\N$ and $n\in\N$ let $(S_d)^n$ denote the $n$-fold
Cartesian product of $S_d$. Define
the mappings $\bf{A},\bf{B}$ from $(S_d)^n$ onto $(S_d)^n$ by
\begin{eqnarray*}
{\bf A}\big( (r^{(1)},r^{(2)} ,\ldots, r^{(n)}) \big) &:=&
\big(\alpha(r^{(1)}),\alpha(r^{(2)}) ,\ldots, \alpha(r^{(n)})\big), \\
{\bf B}\big((r^{(1)},r^{(2)} ,\ldots, r^{(n)})\big) &:=&
\big(\beta(r^{(n)}),\beta(r^{(n-1)}) ,\ldots, \beta(r^{(1)})\big)
\end{eqnarray*}
for $(r^{(1)},r^{(2)} ,\ldots, r^{(n)})\in(S_d)^n$.
Further, let \textbf{id} denote the identity map on $(S_d)^n$. Note that
$\big\{\text{\textbf{id}},\bf{A},\bf{B},\bf{B}\circ \bf{A}\big\}$
together with $\circ$ forms an Abelian group.

We say that
$(Y_k)_{k\in \N}$ is \textit{symmetric in space and time}
iff with respect to each $\vartheta\in\Theta$
\begin{eqnarray*}
\big( Y_{k_1},Y_{k_2},\ldots,Y_{k_n}\big) 
&\stackrel{ \P_{\vartheta} }{=}& \big( -Y_{k_1},-Y_{k_2},\ldots,-Y_{k_n}\big) \\
\text{and} \ \ \ \
\big( Y_{k_1},Y_{k_2},\ldots,Y_{k_n}\big) 
&\stackrel{ \P_{\vartheta} }{=}&
\big( Y_{k_n},Y_{k_{n-1}},\ldots,Y_{k_1}\big) 
\end{eqnarray*}
for all $k_1,k_2,\ldots,k_n\in \N$ with $n\in\N$. 

\begin{lem}\label{revrot}
If  $(X_k)_{k\in\N_0}$ is pairwise distinct and $(Y_k)_{k\in\N}$ 
is symmetric in space and time, then
the distribution of
${\bf\Pi}=(\Pi_d(0),\Pi_d(1),\ldots,\Pi_d(n-1))$ is invariant under 
the group of transformations \rm $\big\{ \{ \text{\textbf{id}}, \bf{A}, \bf{B}, \bf{B}\circ \bf{A} \},\circ\big\}$ \it on $(S_d)^n$
with respect to $\P_H$ for each $\vartheta\in\Theta$, i.e.
\begin{eqnarray*}
{\bf\Pi}\  \stackrel{\P_{\vartheta} }{=}\ {\bf A}({\bf\Pi})\  \stackrel{\P_{\vartheta} }{=}\
{\bf B}({\bf\Pi})\ \stackrel{\P_{\vartheta} }{=}\  {\bf B}\circ {\bf A}({\bf\Pi}).
\end{eqnarray*}
\end{lem}
\begin{proof}
If $(X_k)_{k\in\N_0}$ is pairwise distinct
then, according to the definition of $\widetilde{\pi}$,
\begin{eqnarray*}
&& \big( \widetilde{\pi}((Y_1,Y_2,\ldots,Y_d)), \ldots, \widetilde{\pi}((Y_n,Y_{n+1},\ldots,Y_{d+n-1})) \big)  \\
&& \hspace{7mm} = \ {\bf A}\big(\big( \widetilde{\pi}((-Y_1,-Y_2,\ldots,-Y_d)), \ldots, \widetilde{\pi}((-Y_n,-Y_{n+1},\ldots,-Y_{d+n-1})) \big)\big)  \\
&& \hspace{7mm} = \ {\bf B}\big(\big( \widetilde{\pi}((Y_{d+n-1},\ldots,Y_{n+1},Y_n)), \ldots, \widetilde{\pi}((Y_d,\ldots,Y_{2},Y_1)) \big)\big)  \\
&& \hspace{7mm} = \ {\bf B\circ A}\big(\big( \widetilde{\pi}((-Y_{d+n-1},\ldots,-Y_{n+1},-Y_n)), \ldots, \widetilde{\pi}((-Y_d,\ldots,-Y_{2},-Y_1)) \big)\big) 
\end{eqnarray*}
$\P_{\vartheta}$-a.s. with respect to each $\vartheta\in\Theta$ (see (\ref{patincr})). 
Further, symmetry in space and time of $(Y_k)_{k\in\N}$ implies
that $ (Y_1,Y_2,\ldots,Y_{d+n-1}), 
(-Y_1,-Y_2,\ldots,-Y_{d+n-1}), (Y_{d+n-1},\ldots,Y_2,Y_1)$ and $(-Y_{d+n-1},\ldots,-Y_2,-Y_1) $ 
have the same distribution with respect to $\P_{\vartheta}$ for each $\vartheta\in\Theta$,
and hence the statement follows.
\end{proof}

\paragraph{A Rao-Blackwellization.} For fixed $r\in S_d$ with $d\in\N$ consider the functional $p_r(\cdot)$
defined by $p_r(\vartheta) := \P_{\vartheta}(\Pi_d(0)=r)$ for $\vartheta\in\Theta$. Obviously, 
if $(Y_k)_{k\in\N}$ is stationary then the statistic
\begin{eqnarray}\label{pr}
\widehat{p_{r,}}_n \ = \ \widehat{p_{r,}}_n ({\bf\Pi})  &:=& \frac{1}{n}\sum_{k=0}^{n-1} \1_{ \{ \Pi_d(k) = r \}  } 
\end{eqnarray}
of ${\bf\Pi}=(\Pi_d(0),\Pi_d(1),\ldots,\Pi_d(n-1)) $
is an unbiased estimate of $p_r(\cdot)$,
that is $\E_\vartheta\big( \widehat{p_{r,}}_n \big)=p_r(\vartheta)$ 
for all $\vartheta\in\Theta$ and $n\in\N$.
If, additionally, $(X_k)_{k\in\N_0}$ is pairwise distinct and $(Y_k)_{k\in\N}$
is symmetric in space and time then, according to Lemma
\ref{revrot}, the statistic
\begin{eqnarray}\label{proverline}
\widehat{p_{\overline{r},}}_n &=& 
\widehat{p_{\overline{r},}}_n ({\bf\Pi}) \
:= \
 \frac{1}{4}\Big( \widehat{p_{r,}}_n ({\bf\Pi})  +  \widehat{p_{r,}}_n ({\bf A}({\bf\Pi}))  
+  \widehat{p_{r,}}_n ({\bf B}({\bf\Pi})) +  \widehat{p_{r,}}_n ({\bf B}\circ {\bf A}({\bf\Pi)})      \Big)
\end{eqnarray}
of ${\bf\Pi}$ is a Rao-Blackwellization of $\widehat{p_{r,}}_n$ (see Theorem 3.2.1 in \cite{Pfanzagl94}). This proves the following Theorem.

\begin{theo}\label{RaoBlackwell}
Let $r\in \bigcup_{d\in\N} S_d$. If 
$(X_k)_{k\in\N_0}$ is pairwise distinct and 
$(Y_k)_{k\in\N}$ is stationary and symmetric in space and time
then the estimate
$\widehat{p_{\overline{r},}}_n$ of $p_r(\cdot)$ is unbiased and
more concentrated in convex order than 
$ \widehat{p_{r,}}_n$, that is
\begin{eqnarray}\label{RaoBlackwellconcl}
\E_{\vartheta}(\varphi(
\widehat{p_{\overline{r},}}_n, p_r(\vartheta)  )) &\leq& \E_{\vartheta}(
\varphi(\widehat{p_{r,}}_n , p_r(\vartheta)  ))
\end{eqnarray}
for all $\vartheta\in\Theta$ with respect to each function
$\varphi: [0,1]\times [0,1]\rightarrow [0,\infty[$ with
$\varphi(p,p)=0$ and $\varphi(\cdot,p)$ being convex for every
$p\in [0,1]$.
\end{theo}

According to the strictness of the
Jensen inequality for strictly convex functions, if 
$(X_k)_{k\in\N_0}$ is pairwise distinct and 
$(Y_k)_{k\in\N}$ is stationary and symmetric in space and time,
one has
strict inequality in (\ref{RaoBlackwellconcl})
whenever $\varphi(\cdot,p)$ is strictly convex for every
$p\in [0,1]$ and, additionally,
$\P_{\vartheta}\big( \widehat{p_{\overline{r},}}_n \neq \widehat{p_{r,}}_n \big)
>0$ (see \cite{Pfanzagl94}, Theorem 3.2.1).
Since
$(\cdot-p)^2$ is strictly convex for every
$p\in [0,1]$, in particular
\begin{eqnarray*}
\V_\vartheta
( \widehat{p_{\overline{r},}}_n ) &<& \V_\vartheta ( \widehat{p_{r,}}_n )
\end{eqnarray*}
in this case. Note that for $r\in S_d$ with $d=1$, the estimator 
$\widehat{p_{\overline{r},}}_n$ simply estimates the
constant functional $p_r(\cdot)=\frac{1}{2}$.

\paragraph{Ergodicity.} 
Consider the measurable space
$(\Omega',\cA'):=(\R^{\N},\cB(\R^{\N}))$ of infinite sequences
of real numbers and let the mapping $T:\Omega'\to\Omega'$ be defined
by $T(\omega')=(\omega_2',\omega_3',\ldots)$ for
$\omega'=(\omega_1',\omega_2',\ldots)\in\Omega'$.
The process
$(Y_k)_{k\in \N}$ is called \textit{ergodic} iff 
$(Y_k)_{k\in \N}$ is stationary and, additionally, 
for every $A\in\cA'$
such that $T^{-1}(A)=A$ one has $\P_{\vartheta}((Y_k)_{k\in\N}\in A)=0$ or 
$\P_{\vartheta}((Y_k)_{k\in\N}\in A)=1$ for each $\vartheta\in\Theta$.

As the next Lemma shows, if $(Y_k)_{k\in\N}$ is
ergodic then the estimators of continuous and bounded functionals
of ordinal pattern probabilities are strongly consistent and asymptotically unbiased.

\begin{lem}\label{strongconsist}
Let $r\in \bigcup_{d\in\N} S_d$. If $(Y_k)_{k\in\N}$ is
ergodic and $h:[0,1]\to\R$ is continuous then 
\begin{eqnarray*}
\lim_{n\to\infty}h\big(\widehat{p_{\overline{r},}}_n\big) &=& h\big(p_r(\vartheta)\big)
\end{eqnarray*}
$\P_{\vartheta}$-a.s. for all $\vartheta\in\Theta$. If 
$h$ is continuous and bounded,
then with respect to each $\vartheta\in\Theta$
\begin{eqnarray*}
\lim_{n\to\infty} \E_{\vartheta}\Big( h\big(\widehat{p_{\overline{r},}}_n\big) 
\Big) &=& h\big(p_r(\vartheta)\big).
\end{eqnarray*}
\end{lem}
\begin{proof}
For fixed $\vartheta\in\Theta$ let the probability measure $\mu$ on
$(\Omega',\cA')=(\R^{\N},\cB(\R^{\N}))$ be defined by
$\mu(A):=\P_{\vartheta}((Y_k)_{k\in\N}\in A)$ for $A\in\cA'$.
Further, let the mapping $T:\Omega'\to\Omega'$ be given as above
and define
$f:\Omega'\to\R$ by
\begin{eqnarray*}
f(\omega') &:=& \left\{ \begin{array}{rl}
1&\mbox{for } \widetilde{\pi}((\omega_1',\omega_2',\ldots,\omega_d'))\in \overline{r}\\
0&\mbox{else}
\end{array}\right .
\end{eqnarray*}
for $\omega'=(\omega_1',\omega_2',\ldots)\in\Omega'$. 
Obviously, $f$ is Borel-measurable and  $\int_{\Omega'}|f|d\mu<\infty$,
hence, according to Birkhoff's Theorem, if $(Y_k)_{k\in\N}$ is ergodic then
\begin{eqnarray*}
\lim_{n\to\infty}\widehat{p_{\overline{r},}}_n \ = \ \lim_{n\to\infty}\frac{1}{n}\sum_{j=0}^{n-1}f(T^j((Y_k)_{k\in\N}))
&=& \E_{\vartheta}(f((Y_k)_{k\in\N})) \ = \
p_r(\vartheta)
\end{eqnarray*}
$\P_{\vartheta}$-a.s. (see \cite{Cornfeld82}, Theorem 1.2.1)
and the first statement follows since $h$ is continuous.
If additionally $h$ is bounded then
\begin{eqnarray*}
\lim_{n\to\infty} \E_{\vartheta}\Big( h\big(\widehat{p_{\overline{r},}}_n\big) 
\Big) &=& \E_{\vartheta}\Big( \lim_{n\to\infty} h\big(\widehat{p_{\overline{r},}}_n\big) 
\Big) \ =\ 
h\big(p_r(\vartheta)\big)
\end{eqnarray*}
according to the dominated convergence Theorem.
\end{proof}

\subsection{ Specialization to fBm }\label{fBmspec}

In the following we specialize the considerations of Subsection
\ref{estordp} to equidistant discretizations of fractional Brownian motion
(fBm).
We start with the definition of fBm.

\begin{defi}\label{fBmdefi}
For $H\in\ ]0,1]$ let $(B(t))_{t\in [0,\infty[}$ be a mean-zero
Gaussian process on a probability space $(\Omega,\cA,\P)$ with the
covariance function
\begin{eqnarray*}
\Cov(B(t),B(s))=\frac{1}{2}\big(t^{2H}+s^{2H}-|t-s|^{2H}\big)\, \V(B(1))
\end{eqnarray*}
for $s,t\in [0,\infty[$. Then $(B(t))_{t\in [0,\infty[}$ is called
\textit{fractional Brownian motion (fBm)} with the \textit{Hurst parameter $H$}.
\hfill $\square$ 
\end{defi}

It is well-known that each fBm possesses a
modification with $\P$-a.s. continuous paths (see
\cite{EmbrechtsMaejima02}).
We are only interested
in the distribution of fBm, so let 
for the rest of
the paper $(\P_H)_{H\in ]0,1]}$ be a family of probability measures
on a measurable space $(\Omega,\cA)$, and $(B(t))_{t\in [0,\infty[}$ be a family of real-valued random
variables on $(\Omega,\cA)$ such that $(B(t))_{t\in
[0,\infty[}$ measured with respect to
$\P_H$ is fBm with the Hurst parameter $H$. E.g., $(B(t))_{t\in
[0,\infty[}$ can be defined as the identity on the set
of continuous functions on $[0,\infty[$.

For the \textit{sampling interval length} $\delta>0$ consider
the \textit{equidistant discretization}
\begin{eqnarray*}
(X_k^\delta)_{k\in\N_0} &:=& (B(k\delta))_{k\in\N_0} 
\end{eqnarray*}
of fBm. Further, let $Y_k^\delta:=X_k^\delta-X_{k-1}^\delta$ for $k\in\N$.
Note that fBm is
$H$-self-similar (see \cite{EmbrechtsMaejima02}), that is, for every
$H\in\ ]0,1]$ and for all $a\in\,]0,\infty[$ it holds
\begin{eqnarray*}
(B(at))_{t\in [0,\infty[} &\stackrel{\P_H}{=}&  (a^HB(t))_{t\in
[0,\infty[}
\end{eqnarray*}
with $\stackrel{\P_H}{=}$ denoting equality of all
finite-dimensional distributions with respect to $\P_H$ here.
Since it holds $\widetilde{\pi}((y_1,y_2,\ldots,y_d)) = \widetilde{\pi}((a^H y_1,a^H
y_2,\ldots,a^H y_d))$ for all 
$a>0$, $H\in\ ]0,1]$ and
$(y_1,y_2,\ldots,y_d)\in\R^d$, one has
\begin{eqnarray*}
(\widetilde{\pi}((Y_{k+1}^\delta,\ldots,Y_{k+d}^\delta)))_{k\in\N_0}
&\stackrel{\P_H}{=}&
(\widetilde{\pi}((Y_{k+1}^1,\ldots,Y_{k+d}^1)))_{k\in\N_0}
\end{eqnarray*}
with respect to each $H\in\ ]0,1]$ for all $\delta>0$. This means that 
the distribution of
ordinal patterns for equidistant discretizations of fBm does not
depend on the particular sampling interval length. Therefore, in the
following we only consider the
discretization
$(X_k)_{k\in\N_0} =(X_k^1)_{k\in\N_0}$
with the increment process $(Y_k)_{k\in\N}=(X_k^1-X_{k-1}^1)_{k\in\N}$.

Note that, similarly to the sampling interval length, one can show that
the particular scaling $\V_H(B(1))$ of fBm has no effect on the distribution
of ordinal patterns for equidistant discretizations. Hence,
we always assume the case of \textit{standard fBm} where $\V_H(B(1))=1$
for all $H\in\ ]0,1]$.

According to Definition \ref{fBmdefi}, $\V_H(X_i-X_j)>0$ for all
$i,j\in\N_0$ with $i\neq j$ and $H\in\ ]0,1]$, 
and since $X_i-X_j$ is Gaussian with respect to $\P_H$ for all $H\in\ ]0,1]$, 
the stochastic process $(X_k)_{k\in\N_0}$ is pairwise distinct.
Further, obviously $(Y_k)_{k\in\N}$ is mean-zero
Gaussian with respect to $\P_H$ for each $H\in\ ]0,1]$, and by
Definition \ref{fBmdefi} one obtains
\begin{eqnarray}\label{hdef}
\rho_H(k) \ :=\ \Cov_H(Y_1,Y_{1+k}) &=&
  \frac{1}{2}\big( |k+1|^{2H}-2|k|^{2H}+|k-1|^{2H} \big) 
\end{eqnarray}
for all $H\in\ ]0,1]$ and $k\in\N_0$, that is,
the stochastic process
$(Y_k)_{k\in\N}$ is stationary. Furthermore, for all $n\in\N$
and $k_1,k_2,\ldots,k_n\in\N$
the random vectors $ (Y_{k_1},Y_{k_2},\ldots,Y_{k_n})$, 
$(-Y_{k_1},-Y_{k_2},\ldots,-Y_{k_n})$ and $(Y_{k_n},\ldots,Y_{k_2},Y_{k_1})$
have the same covariance structure, hence $(Y_k)_{k\in\N}$
is space and time symmetric. Consequently,
the conclusion of Theorem \ref{RaoBlackwell} applies to
the estimation of ordinal pattern probabilities for
equidistant discretizations of fBm.

As the following Lemma shows, 
one has strict inequality in (\ref{RaoBlackwellconcl})
whenever $\varphi(\cdot,p)$ is strictly convex for every
$p\in [0,1]$ and $H\in\ ]0,1[$ (compare to the remark after Theorem \ref{RaoBlackwell}).

\begin{lem}
If $H\in\ ]0,1[$ then
$\P_{H}\big( \widehat{p_{\overline{r},}}_n \neq \widehat{p_{r,}}_n \big) > 0$ 
for all $r\in\bigcup_{d\in\N} S_d$ and $n\in\N$.
\end{lem}
\begin{proof}
Let $r\in S_d$ with $d\in\N$ and $n\in\N$ be fixed. It is easy to construct
a sequence of permutations $( r^{(1)}, r^{(2)} ,\ldots, r^{(n)})\in(S_d)^n$
such that the $(n+d-1)$-dimensional Lebesgue measure of
$\bigcap_{k=1}^{n}\big\{\,
(y_1,y_2,\ldots,y_{n+d-1})\in\R^{n+d-1}\ | \ \widetilde{\pi}((y_k,y_{k+1},\ldots,y_{k+d-1}))=r^{(k)} \big\}$
is strictly positive and
\begin{eqnarray*}
\sharp\,\overline{r} \,
\cdot\, \sharp\,\big\{k\in\{0,1,\ldots,n-1\}\,|\, r^{(k)}=r    \big\} 
&\neq& \sharp\,\big\{k\in\{0,1,\ldots,n-1\}\,|\, r^{(k)}\in\overline{r}    \big\}.
\end{eqnarray*}
According to (\ref{hdef}), if $H\in\ ]0,1[$ then $(Y_1,Y_2,\ldots,Y_{n+d-1})$
is non-degenerate Gaussian with respect to $\P_H$ and,
consequently, $\P_H\big( (Y_1,Y_2,\ldots,Y_{n+d-1})\in A \big)>0$
for each Borel-set $A\subset\R^{d+n-1}$ with strictly positive $(n+d-1)$-dimensional Lebesgue measure. Since
$ \widehat{p_{\overline{r},}}_n = \frac{1}{\sharp\,\overline{r}}\sum_{s\in\overline{r}}\widehat{p_{s,}}_n $ 
(see (\ref{pr}) and (\ref{proverline})) the statement follows.
\end{proof}

Note that if $H=1$ then
$\P_{H}\big( \widehat{p_{\overline{r},}}_n \neq \widehat{p_{r,}}_n \big) > 0$
if and only if $r\in\{(0,1,\ldots,d),(d,\ldots,1,0)\}$.

It is well-known that for each $H\in\ ]0,1[$ the spectral
distribution function of $(Y_k)_{k\in\N}$ is absolutely continuous (see \cite{Beran94})
which is a sufficient condition for $(Y_k)_{k\in\N}$ to be ergodic with respect to $\P_H$
(see \cite{Cornfeld82}, Theorem 14.2.1). 
In the case $H=1$
one has $X_k=k X_1$ $\P_H$-a.s. for every $k\in\N_0$ (see \cite{EmbrechtsMaejima02}),
hence
\begin{eqnarray}\label{H=1case}
\widehat{p_{\overline{r}}}_n &=& p_r(H) \ =\
\left\{ \begin{array}{rl}
0&\mbox{if } r\in \bigcup_{d\in\N}\big(S_d\setminus\{(0,1,\ldots,d),(d,\ldots,1,0)\}\big) \vspace{2mm}\\
\frac{1}{2}&\mbox{if } r\in \bigcup_{d\in\N}\{(0,1,\ldots,d),(d,\ldots,1,0)\}
\end{array}\right .
\end{eqnarray}
$\P_H$-a.s. for all $n\in\N$.
Consequently, the conclusions of Lemma \ref{strongconsist} apply to
the estimation of ordinal pattern probabilities for
equidistant discretizations of fBm.

\paragraph{Asymptotic normality.}
Next, we discuss asymptotic normality of the estimators of
ordinal pattern probabilities for equidistant discretizations of fBm. 
We restrict our consideration to the `reasonable' estimators as provided by
Theorem \ref{RaoBlackwell}.
Let $ \stackrel{\P_H}{\longrightarrow}$ denote convergence in distribution with respect
to $\P_H$, and
$N(\mu,\sigma^2)$ be the normal distribution with mean $\mu\in\R$ and
variance $\sigma^2>0$. Given a mean-zero Gaussian vector $ \defn{Z} =(Z_1, Z_2,
\ldots, Z_d)$ on a probability space $(\Omega',\cA',\P)$, and a
Borel-measurable function $f:\R^d\to\R$ such that $\V\big(f(\defn{Z})\big)<\infty$, define the \textit{rank} of $f$ with respect
to $\defn{Z}$ by
\begin{eqnarray*}
\text{rank}(f) \ :=\ \min \big\{
\kappa\in\N: \text{There exists a real polynomial $q:\R^d\to\R$}  \hspace{1.5cm}\\
\text{of degree $\kappa$ with} \ \E\big([f(\defn{Z})-\E(f(\defn{Z}))]\,q(\defn{Z} )\big)\neq 0 \big\},
\end{eqnarray*}
where the minimum of the empty set is infinity. We write $f(k)\sim g(k)$ for mappings $f,g$ from
$\N_0$ onto $\R$ and say that $f$ is \textit{asymptotically equivalent} to $g$
iff $\lim_{k\to\infty}f(k)/g(k)=1$ where $\frac{0}{0}:=1$.
Note that for $\rho_H(k)$ as defined in (\ref{hdef}), one has the asymptotic equivalence 
\begin{eqnarray}\label{hest}
\rho_H(k) &\sim& H(2H-1)k^{2H-2}
\end{eqnarray}
for each $H\in\ ]0,1]$ (see \cite{EmbrechtsMaejima02}).
Further, we write $f(k)=O(g(k))$ iff $\sup_{k\in\N_0}f(k)/g(k)<\infty$,
and $f(k)=o(g(k))$ iff $\lim_{k\to\infty}f(k)/g(k)=0$.

\begin{theo}\label{asnormtheo}
If $H <\frac{3}{4}$ then
\begin{eqnarray*}
\big( \V_H(\widehat{p_{\overline{r},}}_n ) \big)^{-\frac{1}{2}}\,\big( \widehat{p_{\overline{r},}}_n - p_r(H) \big) 
&  \stackrel{\P_H}{\longrightarrow} & 
N(0, 1)
\end{eqnarray*}
for all $r\in S_d$ with $d\in\N\setminus\{1\}$.
\end{theo}
\begin{proof}
Let $r\in S_d$ with $d\in\N\setminus\{1\}$ be fixed and define $f: \R^d\rightarrow \R$ 
by
\begin{eqnarray*}
f((y_1,y_2,\ldots ,y_d)) &:=& \left\{ \begin{array}{rl}
1&\mbox{for
} \widetilde{\pi} ((y_1,y_2,\ldots ,y_d))\in\overline{r} \\
0&\mbox{else}
\end{array}\right .
\end{eqnarray*}
for $(y_1,y_2,\ldots,y_d)\in\R^d$. Write $\defn{Y} = (Y_1,Y_2,\ldots,Y_d)$. We first show that $\text{rank}(f)>1$ with respect to 
$\defn{Y}$ and $\P_H$ for all $H\in\ ]0,1]$. Obviously, $f$ is Borel-measurable,
and $\V_H\big(f( \defn{Y} ) \big)<\infty$ for all $H\in\ ]0,1]$.
Now, let $i\in\{1,2,\ldots,d\}$ be fixed.
Because $(Y_k)_{k\in\N}$ is symmetric in space and time, the random vectors
$(Y_1,Y_2,\ldots,Y_d)$ and $(-Y_1,-Y_2,\ldots,-Y_d)$ 
have the same distribution with respect to $\P_H$ for all $H\in\ ]0,1]$,
and hence, according to the definition of $\alpha$ (see (\ref{alphabeta}))
\begin{eqnarray}\label{proofasnorm}
&& \E_H\big( \1_{ \{ \widetilde{\pi}(\defn{Y})=s\} }  Y_i  \big) 
+ \E_H\big( \1_{ \{ \widetilde{\pi}(\defn{Y})= \alpha(s)\} }  Y_i  \big) \nonumber\\
&& \hspace{20mm} = \
\E_H\big( \1_{ \{ \widetilde{\pi}(\defn{Y})=s\} }  Y_i  \big) 
- \E_H\big( \1_{ \{ \widetilde{\pi}(-\defn{Y})= s\} } (- Y_i)  \big) 
\ = \ 0
\end{eqnarray}
for each $s\in S_d$.
Since $Y_i$ is mean-zero Gaussian with respect to $\P_H$ for all $H\in\ ]0,1]$, 
one has $\E_H(f(\defn{Y}))\,\E_H(Y_i)=0$,
and this together with (\ref{proofasnorm}) yields in case $\sharp\,\overline{r}=2$
\begin{eqnarray*}
\E_H\big([f(\defn{Y})-\E_H(f(\defn{Y}))]\,Y_i\big)\, = \,
\E_H\big( \1_{ \{ \widetilde{\pi}(\defn{Y})=r\} } Y_i \big)  + 
\E_H\big( \1_{ \{ \widetilde{\pi}(\defn{Y})= \alpha(r)\} }
Y_i \big) \, = \, 0
\end{eqnarray*}
with respect to each $H\in\ ]0,1]$, and in case $\sharp\,\overline{r}=4$
\begin{eqnarray*}
 \E_H\big([f(\defn{Y})-\E_H(f(\defn{Y}))]\,Y_i\big) & = &
\E_H\big( \1_{ \{ \widetilde{\pi}(\defn{Y})=r\} } Y_i \big)  + 
\E_H\big( \1_{ \{ \widetilde{\pi}(\defn{Y})= \alpha(r)\} }
Y_i \big)\\ &+&\ 
\E_H\big( \1_{ \{ \widetilde{\pi}(\defn{Y})=\beta(r)\} } Y_i \big)  + 
\E_H\big( \1_{ \{ \widetilde{\pi}(\defn{Y})=\beta\circ\alpha(r)\} }
Y_i \big)
\, = \, 0
\end{eqnarray*}
with respect to each $H\in\ ]0,1]$. Consequently,
$\text{rank}(f)>1$ with respect to $\defn{Y}$ and $\P_H$ for all $H\in\ ]0,1]$. According to (\ref{hest}), $|\rho_H(k)|^{\,\text{rank}(f)} = O\big( k^{4H-4} \big)$
for each $H\in{]0,1]}$, thus 
\begin{eqnarray*}
\sum_{k=0}^\infty
|\rho_H(k)|^{\,\text{rank}(f)} 
\ < \
 \infty
\end{eqnarray*}
for $H < \frac{3}{4}$. Since,
by definition one has $\frac{1}{n}\sum_{k=0}^{n-1}f((Y_{k+1},Y_{k+2}\ldots,Y_{k+d}))=\widehat{p_{\overline{r},}}_n$ 
for every $n\in\N$, 
the statement follows from Theorem 4 of Arcones \cite{Arcones94}.
\end{proof}

Note that Theorem \ref{asnormtheo} can also be proven by the Central Limit
Theorem for non-instanta\-neous filters of a stationary Gaussian
process given by Ho and Sun \cite{HoSun87}.

We leave it as an open question whether $H<\frac{3}{4}$ is also a necessary condition
for $\widehat{p_{\overline{r},}}_n$ to be asymptotically normally distributed.
Indeed, at least for $r\in S_d$ 
with $d=2$ simulations suggest that 
$\widehat{p_{\overline{r},}}_n$
is not asymptotically normally distributed with respect to $\P_H$ if $H\geq\frac{3}{4}$  (see Figure \ref{fig3} below).

\section{ Estimating the probability of a change }\label{var}

\subsection{Ordinal patterns of order \textit{d}=2}

The results of this Section are mainly based on the analysis of 
normal orthant probabilities as given by the following definition.

\begin{defi}\label{Ortprobdef}
Let $n\in\N$ be fixed.
For a non-singular strictly positive definite and symmetric matrix
$\Sigma\in\R^{n\times n}$ let
$\phi(\Sigma,\cdot)$ denote the Lebesgue density of the $n$-dimensional normal
distribution with zero means and the covariance matrix $\Sigma$, that is
\begin{eqnarray*}
\phi(\Sigma, \defn{x}) &=&
\big((2\pi)^n|\Sigma|\big)^{-\frac{1}{2}}\,
\exp\{-\frac{1}{2}\defn{x}^T\Sigma^{-1}\defn{x}
\}
\end{eqnarray*}
for $\defn{x}\in\R^n$. We call
\begin{eqnarray}\label{Phidef}
\Phi(\Sigma) &:=& \int_{{[0,\infty [}^n} \phi(\Sigma,\defn{x})
\,d\defn{x},
\end{eqnarray}
the $n$-\textit{dimensional normal orthant probability} with respect to $\Sigma$.
\hfill $\square$ 
\end{defi}

Clearly, if $(Z_1,Z_2,\ldots,Z_n)$ is a non-degenerate mean-zero Gaussian random vector
on a probability space $(\Omega',\cA',\P)$, then
\begin{eqnarray*}
\Phi\big((\Cov(Z_i,Z_j))_{i,j=1}^n\big) &=& \P(Z_1>0, Z_2>0, \ldots, Z_n>0).
\end{eqnarray*}
The following result is well-known (see \cite{Plackett54}).

\begin{lem}\label{Orthantlemma}
If $(Z_1,Z_2,\ldots,Z_n)$ is a non-degenerate mean-zero Gaussian random vector
on a probability space $(\Omega',\cA',\P)$ such that $\V(Z_k)>0$ for all $k\in\N$, then
\begin{eqnarray*}
\P(Z_1>0, Z_2>0) &=& \frac{1}{4}+\frac{1}{2\pi}\arcsin \rho_{12} , \\
\P(Z_1>0, Z_2>0, Z_3>0) &=& 
\frac{1}{8}+\frac{1}{4\pi}\arcsin \rho_{12} +\frac{1}{4\pi}\arcsin \rho_{13} \\
&& \hspace{22mm}+\frac{1}{4\pi}\arcsin \rho_{23}
\end{eqnarray*}
where $\rho_{ij}=\Corr(Z_i,Z_j)$ for $i,j\in\{1,2,3\}$.
\end{lem}

Note that, in general, no closed-form expressions are available for normal orthant probabilities
of dimension $n\geq 4$. 

Same as in Subsection \ref{fBmspec}, let $(X_k)_{k\in\N_0}$
be an equidistant discretization of fBm with the increment process $(Y_k)_{k\in\N}=(X_k-X_{k-1})_{k\in\N}$.
The following statement is due to Bandt and Shiha \cite{BandtShiha05}.
We refer to parts of the proof below, and thus include it here.

\begin{cor}\label{ordpattprob}
For $H\in\ {]0,1]}$ one has \rm
\begin{eqnarray*}
p_{r}(H) &=& \left\{
\begin{array}{rl}
\frac{1}{\pi}\arcsin 2^{H-1} &\mbox{if }
r\in\{(0,1,2),(2,1,0) \} \vspace{2mm}
\\
\frac{1}{4}-\frac{1}{2\pi}\arcsin 2^{H-1} &\mbox{if }
r\in\{(1,0,2),(1,2,0),(0,2,1),(2,0,1) \}
\end{array}\right . .
\end{eqnarray*}
\end{cor}
\begin{proof}
We first show the statement for $r=(0,1,2)$.
According to (\ref{H=1case}), if $H=1$ then $p_{(0,1,2)}(H)=\frac{1}{2}$.
By formula (\ref{hdef}) one has $\V_H(Y_k)=1$
for each $H\in\ ]0,1]$ and $k\in\N_0$, and hence
hence $\Cov_H(Y_i,Y_j)=\Corr_H(Y_i,Y_j)$ for each $H\in{]0,1]}$
and $i,j\in\N$. Consequently, by
Lemma \ref{Orthantlemma} one obtains
\begin{eqnarray}\label{proofd=2}
p_{(0,1,2)}(H)  & =& \P_H( X_2> X_1> X_0 ) \nonumber\\
&=& \P_H( Y_1> 0, Y_2> 0 ) 
\ =\ \frac{1}{4}+\frac{1}{2\pi}\arcsin \rho_H(1)
\end{eqnarray}
for each $H\in\ {]0,1[}$ with $\rho_H(1)$ as given in (\ref{hdef}). 
Since $\arcsin\rho=2\arcsin\sqrt{(1+\rho)/2}-\frac{\pi}{2}$ for each
$\rho\in[-1,1]$, one gets
$p_{(0,1,2)}(H)=\frac{1}{\pi}\arcsin 2^{H-1}$. Now, by Lemma \ref{revrot} it holds
\begin{eqnarray*}
p_{(0,1,2)}(H)=p_{(2,1,0)}(H) \ \ \ \
\text{and} \ \ \ \
p_{(1,0,2)}(H) = p_{(1,2,0)}(H) = p_{(0,2,1)}(H)  = p_{(2,0,1)}(H)
\end{eqnarray*}
for each $H\in\ ]0,1]$, and hence the statement follows.
\end{proof}

According to Corollary \ref{ordpattprob},
the probabilities of ordinal
patterns of order $d=2$ are all monotonically related to the Hurst
parameter. In particular, consider the \textit{indicator for a change between up to down} 
at time $k$ defined
by
\begin{eqnarray*}
C_k &:=& \1_{\{ X_k \geq X_{k+1} < X_{k+2}\}} +  \1_{\{X_{k} < X_{k+1} \geq X_{k+2}\}}
\end{eqnarray*}
for $k\in\N_0$ and the \textit{probability of a change} 
$c(\cdot)$ given by $c(H) := \P_H(C_{0}=1)$
for $H\in{]0,1]}$. 
Because
$C_k = \1_{\{ \Pi_2(k) = (1,0,2) \}} + \1_{\{ \Pi_2(k) = (1,2,0) \}} 
+ \1_{\{ \Pi_2(k) = (0,2,1) \}} + \1_{\{ \Pi_2(k) = (2,0,1) \}}$
for each $k\in\N_0$, one has
\begin{eqnarray}\label{zetaH}
c(H) &=& p_{(1,0,2)}(H) + p_{(1,2,0)}(H) + p_{(0,2,1)}(H)  + p_{(2,0,1)}(H) \nonumber\\
&=&
1 -\frac{2}{\pi}\arcsin 2^{H-1}
\end{eqnarray}
for $H\in{]0,1]}$. Hence, the probability of a change
and the Hurst parameter $H$ are
monotonically related: the larger $H$, the smaller $c(H)$. In particular,
$c(H)$ tends to $\frac{2}{3}$ as $H$ tends to $0$, and $c(H)=0$ for $H=1$. 

Now, consider the \textit{sample relative frequency of changes} $\widehat{c}_n$
given by
$\widehat{c}_n := \frac{1}{n}\sum_{k=0}^{n-1}C_{k}$ for $n\in\N$.
Since obviously
\begin{eqnarray}\label{estimatoreq}
\widehat{p_{\overline{r},}}_n = \left\{
\begin{array}{rl}
\frac{1}{2}\big( 1 - \widehat{c}_n  \big) &\mbox{if }
r\in\{(0,1,2),(2,1,0) \} \vspace{2mm}
\\
\frac{1}{4}\, \widehat{c}_n &\mbox{if }
r\in\{(1,0,2),(1,2,0),(0,2,1),(2,0,1) \}
\end{array}\right . ,
\end{eqnarray}
in the case $d=2$ any reasonable estimator of ordinal
pattern probabilities as provided by
Theorem \ref{RaoBlackwell} is an affine functional
of $\widehat{c}_n$ and hence 
essentially has the same statistical
properties as $\widehat{c}_n$. In particular, according
to Lemma \ref{strongconsist} and Theorem \ref{asnormtheo}, 
respectively, 
$\widehat{c}_n$ is a strongly consistent estimator of $c(\cdot)$,
and asymptotically normally distributed with respect to $\P_H$
for $H<\frac{3}{4}$. 

In the rest of this Section, we investigate the variance of $\widehat{c}_n$.
In Sec.~\ref{conf} we consider the estimator of $H$ obtained by plugging
the estimate $\widehat{c}_n$ of $c(\cdot)$ 
into the monotonic functional relation (\ref{zetaH}).

Notice that, 
Lemma \ref{Orthantlemma}
yields closed-form expressions
also for probabilities
of ordinal patterns of order $d=3$ (see \cite{BandtShiha05}).
According to Theorem \ref{RaoBlackwell},
one obtains eight different reasonable
estimators of ordinal pattern probabilites in this case where not all 
probabilities seem to be monotonically related to the
Hurst parameter.

\subsection{Variance of the sample relative frequency of changes}\label{Varchanges}
Let $\gamma_H$ denote the autocovariance function
of the stochastic process $(C_k)_{k\in\N_0}$ with respect to $\P_H$ for  $H\in\ ]0,1]$,
that is, $\gamma_H(k) = \Cov_H(C_0,C_{k})$
for  $H\in\ ]0,1]$ and $k\in\N_0$. Then $\V_H(\widehat{c}_n)$
can be written as
\begin{eqnarray}\label{varformula}
\V_H(\widehat{c}_n) &=& \frac{1}{n^2} \big( n\gamma_H(0)+
2\sum_{k=1}^{n-1}(n - k)\gamma_H(k) \big).
\end{eqnarray}
for $n\in\N$. Clearly, $\gamma_H(0)=c(H)(1-c(H))$, and, according to Lemma \ref{Orthantlemma}
one has
\begin{eqnarray}\label{gamma1}
\gamma_H(1)
&=& 2\,\P_H(Y_1> 0, -Y_2> 0, Y_3> 0) - c(H)^2 \nonumber\\
&=& \frac{1}{2\pi}\arcsin\rho_H(2)-\frac{1}{\pi^2}\big(\arcsin\rho_H(1)\big)^2
\end{eqnarray}
for $H\in{]0,1]}$ with $\rho_H(k)$ as given in (\ref{hdef}).
For $k\in\N\setminus\{1\}$ one obtains
\begin{eqnarray}\label{gamma2}
\gamma_H(k) &=& 
\Cov_H(1-C_0,1-C_{k}) \nonumber\\\nonumber
&=& \P_H(C_0=0,C_k=0)-\P_H(C_0=0)^2 \\\nonumber
&=& 2\, \P_H(Y_1> 0,Y_2>0,Y_{k+1}>0,Y_{k+2}>0) \\\nonumber
&& \ \ + \ 2\, \P_H(Y_1>0,Y_2>0,-Y_{k+1}> 0,-Y_{k+2}> 0) \\
&& \ \ \ \ \ \ \ \ - \ 4\, \P_H(Y_1>0,Y_2>0)^2
\end{eqnarray}
for $H\in{]0,1]}$ where, in general,
no closed-form expressions are known
for the probabilities $\P_H(Y_1>0,Y_2>0,Y_{k+1}>0,Y_{k+2}>0)$ and
$\P_H(Y_1>0,Y_2>0,-Y_{k+1}> 0,-Y_{k+2}> 0)$, respectively.

Obviously, if $H=1$ then $\widehat{c}_n=0$ $\P_H$-a.s. 
and hence $\V_H(\widehat{c}_n)=0$ for all $n\in\N$ (compare to 
(\ref{estimatoreq}) and (\ref{H=1case})).
Note that, indeed, $\gamma_H(k)=0$ for each $k\in\N_0$ in this case
which can be seen from the fact that $Y_i=Y_j$ $\P_H$-a.s. 
for all $i,j\in\N$ if $H=1$.

In the case $H=\frac{1}{2}$ one obtains $\gamma_H(0)=\frac{1}{4}$ and $\gamma_H(k)=0$ for all $k\in\N$,
and hence, according to formula (\ref{varformula}),
$\V_H(\widehat{c}_n)=\frac{1}{4n}$ for $n\in\N$.
In particular, $\rho_H(k)=0$ for all $k\in\N$ if $H=\frac{1}{2}$ (see (\ref{hdef})),
hence the $Y_1,Y_2,Y_{k+1},Y_{k+2}$ in (\ref{gamma2}) are stochastically independent
and thus
\begin{eqnarray*}
&& \P_H(Y_1> 0,Y_2>0,Y_{k+1}>0,Y_{k+2}>0) \\
&& \hspace{12mm}= \ \P_H(Y_1>0,Y_2>0,-Y_{k+1}> 0,-Y_{k+2}> 0) \ = \
\P_H(Y_1>0,Y_2>0)^2 \ = \ \frac{1}{16}.
\end{eqnarray*}

\paragraph{Numerical evaluation of the covariances.}
Next, we provide a way for the numerical evaluation of $\gamma_H(k)$ 
in the non-trivial case $H\notin\{\frac{1}{2},1\}$ and 
$k\in\N\setminus\{1\}$. 
Let $\cR$ denote the set of
$\defn{r}=(r_1,r_2,r_3,r_4)\in[-1,1]^4$ such that the matrix
\begin{eqnarray}\label{Sigma}
\Sigma(\defn{r}) &:=&\left(
\begin{array}{cccc}
1 & r_{1} & r_{2} & r_{3} \\
r_{1} & 1 & r_{4} & r_{2} \\
r_{2} & r_{4} & 1 & r_{1} \\
r_{3} & r_{2} & r_{1} & 1
\end{array}
 \right)
\end{eqnarray}
is strictly positive definite. Note that, if $\defn{r},\defn{s}\in\cR$ 
then $\defn{x}^T\Sigma(h\cdot\defn{r}+(1-h)\cdot\defn{s})\,\defn{x} = 
h\,\defn{x}^T\Sigma(\defn{r})\,\defn{x}+
(1-h)\,\defn{x}^T\Sigma(\defn{s})\,\defn{x} > 0$
for each $\defn{x}\in\R^4\setminus\{\defn{0}\}$ and $h\in[0,1]$, hence $\cR$ is convex.
Define
\begin{eqnarray}\label{rho}
\rho((k,H,h)) &:=& (\rho_H(1),h\cdot\rho_H(k),h\cdot\rho_H(k+1),h\cdot\rho_H(k-1)) 
\end{eqnarray}
for $k\in\N\setminus\{1\}$, $H\in\ ]0,1[$ and $h\in[-1,1]$. 
Obviously, $\Sigma\big(\rho((k,H,1))\big)$ 
is the covariance matrix of $(Y_1,Y_2,Y_{k+1},Y_{k+2})$,
and $\Sigma\big(\rho((k,H,-1))\big)$ is the covariance matrix of $(Y_1,Y_2,-Y_{k+1},-Y_{k+2})$, both with respect to $\P_H$, and since $(Y_1,Y_2,Y_{k+1},Y_{k+2})$
and $(Y_1,Y_2,-Y_{k+1},-Y_{k+2})$ are non-degenerate Gaussian with respect to $\P_H$
if $H\in\ ]0,1[$, one has $\rho((k,H,1))\in\cR$
and $\rho((k,H,-1))\in\cR$ for each $k\in\N\setminus\{1\}$ and $H\in\ ]0,1[$. From
the convexity of
$\cR$ it follows that $\rho((k,H,h))\in\cR$ for each $k\in\N\setminus\{1\}$, $H\in\ ]0,1[$ and $h\in[-1,1]$,
consequently,
\begin{eqnarray*}
v((k,H,h)) &:=& (\Phi\circ\Sigma)\big( \rho((k,H,h)) \big)
\end{eqnarray*}
is well-defined for all $ k\in\N\setminus\{1\}$, $H\in\ ]0,1[$ and $h\in[-1,1]$
(compare to Definition \ref{Ortprobdef}). 
Note that
\begin{eqnarray*}
v((k,H,1)) &=&\P_H(Y_1 > 0,Y_2 > 0,Y_{k+1} > 0,Y_{k+2} > 0),\\
v((k,H,-1)) &=& \P_H(Y_1 > 0,Y_2 >0,-Y_{k+1} > 0,-Y_{k+2} > 0), \\
v((k,H,0)) &=& \P_H(Y_1 > 0,Y_2 > 0)^2 
\end{eqnarray*}
for all $ k\in\N\setminus\{1\}$ and $H\in\ {]0,1[}$. Hence, by inserting 
the expressions on the left side
into (\ref{gamma2}) one obtains
\begin{eqnarray}\label{numform}
\gamma_H(k) &=&
2\, \int_0^1 \Big(\, \frac{\partial v}{\partial h}((k,H,x)) 
- \frac{\partial v }{\partial h}((k,H,-x)) \,\Big)\, dx
\end{eqnarray}
for $ k\in\N\setminus\{1\}$ and $H\in\ ]0,1[$. 
The partial derivative of $v$ with respect to $h$ is given by
\begin{eqnarray*}
 \frac{\partial v}{\partial h}((k,H,x)) & = &
\rho_H(k) \frac{\partial (\Phi\circ\Sigma) }{\partial r_2}\big(\rho((k,H,x))\big) 
\, + \,
\rho_H(k+1) \frac{\partial (\Phi\circ\Sigma)}{\partial r_3} \big(\rho((k,H,x))\big) \\
&& \hspace{20mm}
+\, \rho_H(k-1) \frac{\partial (\Phi\circ\Sigma)}{\partial r_4}\big(\rho((k,H,x))\big)
\end{eqnarray*}
for $k\in \N\setminus\{1\}$, $H\in\ ]0,1[$ and $x\in[-1,1]$. 
By the reduction formula for normal orthant
probabilities given by Plackett (see \cite{Plackett54}), one gets the first partial
derivatives of $(\Phi\circ\Sigma)$ with respect to $r_2,r_3,r_4$, namely 
\begin{eqnarray}\label{firstpartder}
\frac{\partial (\Phi\circ\Sigma)}{\partial r_2}(\defn{s}) &=& \frac{1}{\pi(1-s_2^2)^{\frac{1}{2}}}
\left( \frac{1}{4}-\frac{1}{2\pi}\arcsin \frac{ |\Sigma(\defn{s})_{13}| }{(|\Sigma(\defn{s})_{11}| |\Sigma(\defn{s})_{22}|)^{\frac{1}{2}} }\right), \nonumber\\
\frac{\partial (\Phi\circ\Sigma)}{\partial r_3}(\defn{s}) &=& \frac{1}{2\pi(1-s_3^2)^{\frac{1}{2}}}
\left( \frac{1}{4}+\frac{1}{2\pi}\arcsin \frac{ |\Sigma(\defn{s})_{23}| }{ |\Sigma(\defn{s})_{22}| }\right), \nonumber\\ 
\frac{\partial ( \Phi\circ\Sigma )}{\partial r_4}(\defn{s}) &=& \frac{1}{2\pi(1-s_4^2)^{\frac{1}{2}}}
\left( \frac{1}{4}+\frac{1}{2\pi}\arcsin \frac{ |\Sigma(\defn{s})_{14}| }{ |\Sigma(\defn{s})_{11}| }\right)
\end{eqnarray}
for $\defn{s}=(s_1,s_2,s_3,s_4)\in\cR$, 
where $\Sigma(\defn{s})_{ij}$ denotes the matrix obtained from $\Sigma(\defn{s})$ by deleting the $i$-th row and $j$-th column of $\Sigma(\defn{s})$.

Formula (\ref{numform}) together with
formula (\ref{varformula})
allows to compute numerical values of $\V_H(\widehat{c}_n)$ 
to any desired precision.
See \cite{Cheng69} for details on the evaluation of the integral in (\ref{numform}).
The following Lemma will be needed below.

\begin{lem}\label{gammacontlemma}
For each $k\in\N_0$ the mapping $H\mapsto \gamma_H(k)$ is continuous on $]0,1[$.
\end{lem}
\begin{proof}
For $k=0$ the statement is valid since $H\mapsto c(H)(1-c(H))$ is continuous on $]0,1]$.
Note that $H\mapsto \rho_H(k)$ is continuous on $]0,1]$
for each $k\in\N$ (see (\ref{hdef})). Consequently, the statement is valid
for $k=1$ (compare to (\ref{gamma1})) and, furthermore, 
the mapping $H\mapsto v((k,H,h))$ is continuous on $]0,1[$
for each $k\in\N\setminus\{1\}$ and $h\in[-1,1]$.
Since $\gamma_H(k)=2v((k,H,1))+2v((k,H,-1))-4v((k,H,0))$,
the statement follows.
\end{proof}

\subsection{Limit behaviour}
The next Theorem 
establishes
an asymptotically equivalent expression for $\gamma_H(k)$.

\begin{theo}\label{ascov}
For each $H\in\ ]0,1]$ one has
\begin{eqnarray*}
\gamma_H(k) &\sim&  \frac{2(1-\rho_H(1))}{\pi^2(1+\rho_H(1))} H^2(2H-1)^2k^{4H-4}
\end{eqnarray*}
as $k$ tends to $\infty$.
\end{theo}
\begin{proof}
First note, for $H=\frac{1}{2}$ and $H=1$ the statement is true since 
both sides evaluate to $0$ for each $k\in\N$ (in particular $\rho_H(1)=1$ for $H=1$).
Let the mapping $(\Phi\circ\Sigma): \cR\to[0,1]$ be given as in (\ref{Phidef})
and (\ref{Sigma}), respectively. 
For each $\defn{s}=(s_1,s_2,s_3,s_4)\in\cR$
Taylor's Theorem asserts  
the existence of some $h\in [0, 1]$ such that
\begin{eqnarray*}
(\Phi\circ\Sigma)(\defn{s}) &=& (\Phi\circ\Sigma)((s_1,0,0,0))
+ \sum_{i=2}^4 s_i \frac{\partial ( \Phi\circ\Sigma )}{\partial r_i}((s_1,0,0,0)) \nonumber\\
&& \hspace{27mm}
+ \hspace{2mm} \frac{1}{2}\sum_{i,j=2}^4 s_is_j \frac{\partial^2 ( \Phi\circ\Sigma )}{\partial r_i \partial r_j}((s_1,0,0,0)) \nonumber \\ 
&& \hspace{22mm}
+ \hspace{2mm} \frac{1}{6}\,\sum_{i,j,l=2}^4 s_is_js_l 
\frac{\partial^3 ( \Phi\circ\Sigma )}{\partial r_i \partial r_j \partial r_l}
((s_1,h\cdot s_2,h\cdot s_3,h\cdot s_4)).
\end{eqnarray*}
Consequently, for each $(s_1,s_2,s_3,s_4)\in\cR$ with $(s_1,-s_2,-s_3,-s_4)\in\cR$
there exist $h_1\in [0, 1]$ and $h_2\in [-1, 0]$ such that
\begin{eqnarray}\label{Taylorexpr}
(\Phi\circ\Sigma)( (s_1,s_2,s_3,s_4) ) \ +\ (\Phi\circ\Sigma)( (s_1,-s_2,-s_3,-s_4) )
\ = \hspace{30mm}&& \nonumber\\
 \hspace{-10mm} 2\,(\Phi\circ\Sigma)((s_1,0,0,0)) \
+ \ \sum_{i,j=2}^4 s_is_j \frac{\partial^2 ( \Phi\circ\Sigma )}{\partial r_i \partial r_j}((s_1,0,0,0)) && \nonumber \\ 
+ \hspace{2mm} \frac{1}{3}\,\sum_{i,j,l=2}^4 s_is_js_l 
\frac{\partial^3 ( \Phi\circ\Sigma )}{\partial r_i \partial r_j \partial r_l}
((s_1,h_1\cdot s_2,h_1\cdot s_3,h_1\cdot s_4)) && \nonumber \\ 
- \hspace{2mm} \frac{1}{3}\,\sum_{i,j,l=2}^4 s_is_js_l 
\frac{\partial^3 ( \Phi\circ\Sigma )}{\partial r_i \partial r_j \partial r_l}
((s_1,h_2\cdot s_2,h_2\cdot s_3,h_2\cdot s_4)) && \hspace{-5mm}.
\end{eqnarray}
Now, according to equation (\ref{gamma2}), for $H\in\ ]0,1[$ and $k\in\N\setminus\{1\}$
one has
\begin{eqnarray}\label{gammacont}
\gamma_H(k) 
&=& 2\, (\Phi\circ\Sigma) \big(( \rho_H(1), \rho_H(k), \rho_H(k+1), \rho_H(k-1) )\big) \nonumber\\ 
&& \ \ + \ 2\, (\Phi\circ\Sigma) \big(( \rho_H(1), -\rho_H(k), -\rho_H(k+1), -\rho_H(k-1) )\big)  \nonumber \\
&& \ \ \ \ \ \ \ \ - \ 4\, (\Phi\circ\Sigma) \big(( \rho_H(1), 0,0,0 )\big).
\end{eqnarray}
Let $ \widetilde{\rho}_H(k) :=H(2H-1)k^{2H-2} $ for $H\in{]0,1]}$ and $k\in\N\setminus\{1\}$.
By formula (\ref{hest}), $\rho_H(k)$ is asymptotically equivalent
to $ \widetilde{\rho}_H(k)$ for each $H\in{]0,1]}$, hence 
for all $i,j,l\in\N$ one has the asymptotic equivalence 
\begin{eqnarray*}
(\rho_H(k))^i\,(\rho_H(k+1))^j\,(\rho_H(k-1))^l &\sim&
(\widetilde{\rho}_H(k))^i\,(\widetilde{\rho}_H(k))^j\,(\widetilde{\rho}_H(k))^l
\end{eqnarray*}
for each $H\in{]0,1]}$. Consequently,
according to formula (\ref{Taylorexpr}) and (\ref{gammacont}), respectively,
for each $H\in\ ]0,1[$ and $k\in\N\setminus\{1\}$ there exist
$h_1\in [0, 1]$ and $h_2\in [-1, 0]$ such that
\begin{eqnarray}\label{ascovproof}
\gamma_H(k) 
&\sim& 2\,\big(\widetilde{\rho}_H(k)\big)^2\, 
\sum_{i,j=2}^4
\frac{\partial^2 ( \Phi\circ\Sigma )}{\partial r_i \partial r_j} \big(( \rho_H(1), 0,0,0 )\big) \
 + \  2\,\big(\widetilde{\rho}_H(k)\big)^3\, R_H(k)
\end{eqnarray}
where
\begin{eqnarray*}
R_H(k) &=& 
 \frac{1}{3}\sum_{i,j,l=2}^4
\frac{\partial^3 ( \Phi\circ\Sigma )}{\partial r_i \partial r_j \partial r_l}
\big( \rho((k,H,h_1)) \big) 
- \frac{1}{3}\sum_{i,j,l=2}^4 
\frac{\partial^3 ( \Phi\circ\Sigma )}{\partial r_i \partial r_j \partial r_l}
\big( \rho((k,H,h_2)) \big)
\end{eqnarray*}
with $\rho((k,H,h))$ as defined in (\ref{rho}).
Now, by
the formulas for the first partial derivatives of 
$ (\Phi\circ\Sigma) $ given in (\ref{firstpartder}) 
one obtains
\begin{eqnarray*}
&& 
\frac{\partial^2 ( \Phi\circ\Sigma )}{\partial^2 r_2} ( \defn{s} ) 
=\frac{1+s_1^2}{2\pi^2(1-s_1^2)}, \ \ \
\frac{\partial^2 ( \Phi\circ\Sigma )}{\partial r_2 \partial r_3} ( \defn{s} )  
=\frac{\partial^2 ( \Phi\circ\Sigma )}{\partial r_2 \partial r_4} ( \defn{s} )
=-\frac{s_1}{2\pi^2(1-s_1^2)}, \\
&&
\frac{\partial^2 ( \Phi\circ\Sigma )}{\partial^2 r_3} ( \defn{s} ) 
=\frac{\partial^2 ( \Phi\circ\Sigma )}{\partial^2 r_4} ( \defn{s} ) 
=\frac{s_1^2}{4\pi^2(1-s_1^2)}, \ \ \
\frac{\partial^2 ( \Phi\circ\Sigma )}{\partial r_3 \partial r_4} ( \defn{s} )  
=\frac{1}{4\pi^2(1-s_1^2)}
\end{eqnarray*}
for $\defn{s}=(s_1,0,0,0)\in\cR$. In particular, one has
$|\Sigma( \defn{s} )_{13}|=|\Sigma(\defn{s})_{23}|=|\Sigma(\defn{s})_{14}|=0$
for $\defn{s}=(s_1,0,0,0)\in\cR$.
Putting these terms with $(s_1,0,0,0)=(\rho_H(1),0,0,0)$ into (\ref{ascovproof}),
one gets
\begin{eqnarray*}
\gamma_H(k) \ 
\sim\ 2\,\big(\widetilde{\rho}_H(k)\big)^2\, 
\frac{4 - 8\rho_H(1) + 4\rho_H(1)^2}{4\pi^2(1 - \rho_H(1)^2)}
\,+ \,  2\,\big(\widetilde{\rho}_H(k)\big)^3\, R_H(k)&&  \\
=\ \frac{2(1 - \rho_H(1))}{\pi^2(1 + \rho_H(1))}\,H^2(2H-1)^2k^{4H-4}\,+ \,  
O(k^{6H-6}\, R_H(k))&&
\end{eqnarray*}
for each $H\in\ ]0,1[$. The final claim is $k^{6H-6}\,R_H(k)=o(k^{4H-4})$ for each $H\in\ {]0,1[}$.
Let $H\in\ {]0,1[}$ be fixed.
Since $\rho((k,H,1))$ and $\rho((k,H,-1))$ both converge to $(\rho_H(1),0,0,0)$
as $k$ tends to $\infty$,
the convex hull of
the subset
\begin{eqnarray*}
\cR_H &:=& 
\Big(\,\bigcup_{k\in\N\setminus\{ 2\}}\big\{\rho((k,H,1)), \rho((k,H,-1))\big\}\,\Big)\cup
\big\{(\rho_H(1),0,0,0) \big\}
\end{eqnarray*}
of $\cR$ is closed in $[-1,1]^4$.
Note that $\bigcup_{k\geq 2}\big\{\rho((k,H,x)): \ x\in[-1,1]  \big\}$
is a subset of the convex hull of $\cR_H$. Consequently,
since $\frac{\partial^3 ( \Phi\circ\Sigma )}{\partial r_i \partial r_j \partial r_l}$
is continuous on $\cR$ for all $i,j,l\in\{2,3,4\}$ (see \cite{Plackett54}),
one has
\begin{eqnarray*}
\sup_{k\in\N\setminus\{ 2\}}R_H(k) < \infty
\end{eqnarray*}
which shows the claim.
\end{proof}

In the case $H\geq\frac{3}{4}$, equation (\ref{varformula}) together with Theorem \ref{ascov} yields
simple closed-form asymptotically equivalent expressions for $\V_H( \widehat{c}_n )$, namely
\begin{eqnarray}\label{varcn}
\V_H( \widehat{c}_n ) &\sim& \left\{ \begin{array}{rl}
d_H\frac{\ln(n)}{n} &\mbox{if } H=\frac{3}{4}\vspace{2mm}
\\
d_H\frac{n^{4H-4}}{4H-3} &\mbox{if } H>\frac{3}{4}
\end{array}\right .
\end{eqnarray}
where
$d_H = \frac{4(1-\rho_H(1))}{\pi^2(1+\rho_H(1))}H^2(2H-1)^2$. As one can see,
$\V_H( \widehat{c}_n )$ tends to $0$ with a rate slower than $n^{-1}$ in this case.
Note that, from the practical viewpoint, the expressions in (\ref{varcn}) are
useful approximations of $\V_H(\widehat{c}_n)$ only if $n$ is excessively large.
We provide a way for getting more accurate practical approximations at the end of Sec.~\ref{conf}. 

If $H<\frac{3}{4}$ then, 
by equation (\ref{varformula}) and Theorem \ref{ascov}, one has
$\V_H( \widehat{c}_n )=O(n^{-1})$.
Let
a family of mappings $(f_n)_{n\in\N}$ from $]0,\frac{3}{4}[$ onto $\R$ 
be defined by 
\begin{eqnarray}\label{deffn}
f_n(H) &:=& n\cdot\V_H(\widehat{c}_n)
\end{eqnarray}
for $H\in {]0,\frac{3}{4}[}$ and $n\in\N$. 
Further, define
$f_\infty(H):= \lim_{n\to\infty} f_n(H)$
for $H\in \ {]0,\frac{3}{4}[}$. 
The following Lemma will be needed in the next section. 

\begin{lem}\label{varconv}
$f_n$ uniformly converges to $f_\infty$ on every compact subset $I\subset\ ]0,\frac{3}{4}[$.
Moreover, $f_\infty$ is continuous and strictly positive on $]0,\frac{3}{4}[$.
\end{lem}
\begin{proof}
First, by Theorem \ref{ascov} one has $\sum_{k=1}^\infty\gamma_H(k)<\infty$ for each 
$H\in \ {]0,\frac{3}{4}[}$. Hence, according to equation (\ref{varformula})
and the dominated convergence Theorem one obtains
\begin{eqnarray}\label{proofvarconv2}
f_\infty(H) - \gamma_H(0)&=& \lim_{n\to\infty} 2\sum_{k=1}^{n-1} \frac{n-k}{n}\,\gamma_H(k)
\ = \  \lim_{n\to\infty} 2\sum_{k=1}^{\infty} \max\big\{0,\frac{n-k}{n}\}\,\gamma_H(k) 
\nonumber\\
&=& 2\sum_{k=1}^{\infty} \lim_{n\to\infty} \max\big\{0,\frac{n-k}{n}\}\,\gamma_H(k) 
\ =\  2\sum_{k=1}^\infty\gamma_H(k)
\end{eqnarray}
and, consequently,
\begin{eqnarray}\label{proofvarconv1}
f_\infty(H) - f_n(H) &=& 
2\,\Big( \,\sum_{k=1}^\infty \gamma_H(k)  - \sum_{k=1}^{n-1} \frac{n-k}{n}\, \gamma_H(k) \Big) \nonumber\\
&=& 2\, \sum_{k=1}^\infty \min\big\{k/n,1\big\}\, \gamma_H(k)
\end{eqnarray} 
for each $H\in \ {]0,\frac{3}{4}[}$. Now, let $I\subset\ ]0,\frac{3}{4}[$ be compact. 
By Lemma \ref{gammacontlemma}, the mapping $H\mapsto \gamma_H(k)$
is continuous on $I$ for every $k\in\N$, hence one can define the family
$(H_k)_{k\in\N}$ of numbers in $I$ to be such that
$\gamma_{H_k}(k) = \max\{\gamma_H(k): \ H\in I\}$ for every $k\in\N$. 
Then one has
\begin{eqnarray}\label{unif}
\big| f_\infty(H) - f_n(H) \big| &\leq& 2\, \sum_{k=1}^\infty \min\big\{k/n,1\big\}\,
\big|\gamma_{H_k}(k) \big|
\end{eqnarray} 
uniformly for all $H\in I$. 
To prove the uniform convergence of $f_n$ to $f_\infty$ on $I$,
it is sufficient to show that the right side of (\ref{unif})
tends to $0$ as $n$ tends to $\infty$. Let $\delta\in{]0,\frac{1}{4}[}$ 
and define $H_{\ast}:= 
\max \big(I\cup\{ \frac{1}{2}+\delta\}\big)$.
By Theorem \ref{ascov} one has
$|\gamma_{H_k}(k) |=O\big(|\gamma_{H_{\ast}}(k) |\big)$,
and hence
\begin{eqnarray*}
\sum_{k=1}^\infty \min\big\{k/n,1\big\}\,
\big|\gamma_{H_k}(k) \big|
\,=\, O\Big(\sum_{k=1}^\infty \min\big\{k/n,1\big\}\,
\big|\gamma_{H_{\ast}}(k) \big|\Big) .
\end{eqnarray*}
With the same arguments as in (\ref{proofvarconv2}) one obtains
\begin{eqnarray*}
\lim_{n\to\infty}\sum_{k=1}^\infty \min\big\{k/n,1\big\}\,
\big|\gamma_{H_{\ast}}(k) \big| &=&
\sum_{k=1}^\infty \lim_{n\to\infty}\min\big\{k/n,1\big\}\,
\big|\gamma_{H_{\ast}}(k) \big| \ = \ 0 
\end{eqnarray*}
showing that the right side of (\ref{unif})
tends to $0$ as $n$ tends to $\infty$.
Now, the continuity of $f_\infty$ follows since
$f_n$ is continuous on $]0,\frac{3}{4}[$ for every
$n\in\N$ which is a direct consequence of formula (\ref{varformula}) together with Lemma \ref{gammacontlemma}.
To establish the strict positiveness of $f_\infty$,
first note that $f_\infty(\frac{1}{2})=\frac{1}{4}$ (see Subsec.~\ref{Varchanges}). 
Further, let $H\in\ ]0,\frac{3}{4}[\setminus\{\frac{1}{2}\}$ 
be fixed. According to Theorem \ref{ascov},
one has asymptotically
$\gamma_H(k)>0$, hence, by equation (\ref{proofvarconv1}), there exists an $n_0\in\N$ such that
$f_\infty(H) - f_{n_0}(H)>0$. Obviously, $f_n(H)\geq 0$ for each $H\in\ ]0,\frac{3}{4}[$ 
and $n\in\N$ which completes the proof.
\end{proof}

\section{Confidence intervals for the ZC estimator}\label{conf} 

According to equation (\ref{zetaH}), the probability of a change $c(H)$
and the Hurst parameter $H$ are
monotonically related. Here we consider the estimator $\widehat{H}_n$
of $H$ obtained by plugging
the estimate $\widehat{c}_n$ of $c(\cdot)$ 
into the monotonic functional relation (\ref{zetaH}).
In order to get finite non-negative estimates, we define
\begin{eqnarray}\label{HfromZeta}
g(x) &:=& \left\{
\begin{array}{rl}
\log_2\sin(\pi(1-x)/2)+1 & \mbox{for
} x\in\big[0, \frac{2}{3} \big[ \vspace{2mm} \\
0&\mbox{for
} x\in\big[\frac{2}{3},1 \big]
\end{array}
\right.
\end{eqnarray}
and set $\widehat{H}_n := g(\widehat{c}_n)$. 
Obviously, the first derivative $g'$ of $g$ exists 
and is non-zero on $]0,\frac{2}{3}[$, namely one has
\begin{eqnarray}\label{gder}
g'(x) &=& -\frac{\pi}{2\ln 2}\frac{ \cos(\pi(1-x)/2 )}{\sin(\pi(1-x)/2 )}  \ < \ 0
\end{eqnarray}
for $x\in\ ]0,\frac{2}{3}[$.
Hence, by
Theorem \ref{asnormtheo} together with Theorem 2.5.2 in \cite{Lehmann99}
one obtains
\begin{eqnarray*}
\big( g'(c(H))\big)^{-1}\,\big( \V_H(\widehat{c}_n) \big)^{-\frac{1}{2}} \,\big( \widehat{H}_n - H \big) 
&  \stackrel{\P_H}{\longrightarrow} & 
N(0,1)
\end{eqnarray*}
for each $H\in\ ]0,\frac{3}{4}[$. Note that $g$ is continuous on $[0,1]$, and thus
$\lim_{n\to\infty}\widehat{H}_n=H$ $\P_H$-a.s. for each $H\in{]0,1]}$ (see Lemma \ref{strongconsist}).
Let $f_n(H)$ be defined as in (\ref{deffn}).
According to Lemma \ref{varconv}, and because the mapping $H\mapsto g'(c(H))$ is continuous on $]0,1[$, one obtains
\begin{eqnarray*}
\lim_{n\to\infty}\,
\frac{ \big( g'(c( \widehat{H}_n    ))\big)^2\,f_n( \widehat{H}_n ) }{\big( g'(c(H))\big)^2\,f_n(H)} &=& 
\frac{ \big( g'(c(H))\big)^2 \,f_\infty( H)}{\big( g'(c(H))\big)^2\,f_\infty(H)} \ =\ 1
\end{eqnarray*}
$\P_H$-a.s. for each $H\in\ ]0,\frac{3}{4}[$. 
By Theorem 2.3.3 in \cite{Lehmann99} it follows that
\begin{eqnarray}\label{Hasnorm}
\big( g'(c( \widehat{H}_n    ))\big)^{-1}\,
\big( \big( f_n(\widehat{H}_n) \big)^{-\frac{1}{2}}\,
\sqrt{n} \,\big( \widehat{H}_n - H \big) 
&  \stackrel{\P_H}{\longrightarrow} & N(0,1)
\end{eqnarray}
for each $H\in\ ]0,\frac{3}{4}[$. Now, let the family of mappings $(s_n)_{n\in\N}$ from
$[0,1]$ on $\R$ be defined by $s_n(1):=0$,
\begin{eqnarray*}
s_n(H) &:=& n^{-1} \big( g'(c(H))\big)^2 \V_H(\widehat{c}_n)
\end{eqnarray*}
for $H\in{]0,1[}$, and $s_n(0):=\lim_{H\to 0}s_n(H)$.
Note that $s_n(0)$ is well-defined for each $n\in\N$ since,
according to Lemma \ref{gammacontlemma} and equation (\ref{gder}),
the mappings $H\mapsto\V_H(\widehat{c}_n)$ 
and $H\mapsto g'(c(H))$, respectively, are
continuous on $]0,1[$. For $n\in\N$ consider the random (confidence)
interval 
\begin{eqnarray}\label{convintform}
K_n &:=& \big[ \,\widehat{H}_n-1.96 \, \big( s_n(\widehat{H}_n) \big)^{\frac{1}{2}}, \,\, 
\widehat{H}_n+1.96 \, \big( s_n(\widehat{H}_n) \big)^{\frac{1}{2}}\, \big]
\,\cap\,\big[0,1\big].
\end{eqnarray}
Let $\Phi$ denote the cumulative distribution function of
$N(0,1)$ here. According to (\ref{Hasnorm}) one obtains
\begin{eqnarray}\label{confidence1}
\P_H\big( H \in K_n  \big) 
& \longrightarrow & 2\,\Phi(1.96)-1 \ = \ 0.9500...
\end{eqnarray}
as $n$ tends to $\infty$. 
We compare this to simulations of
$\P_H\big( H\in K_n \big)$ below (see Table \ref{tab2}).
The computations of $s_n(\widehat{H}_n)$
were carried out using formula (\ref{varapprox}) given below.
See Figure \ref{fig1} (a) for a visualization of
the confidence intervals.

\paragraph{Asymptotic expectation and variance.}
Notice that, simulations suggest that
the actual expectation $\E_H(\widehat{H}_n)$ and
the actual variance 
$\V_H( \widehat{H}_n )$ 
of $\widehat{H}_n$
asymptotically behave like the expectation and the variance
of the asymptotic distribution of $\widehat{H}_n$ obtained by the Delta method (see \cite{Lehmann99}),
namely, 
\begin{eqnarray}\label{asbias}
\E_H(\widehat{H}_n) &\sim&  H-\frac{1}{2}\,g''(c(H))\V_H( \widehat{c}_n  )
\end{eqnarray}
as well as 
\begin{eqnarray}\label{asvar}
\V_H(\widehat{H}_n) &\sim& g'(c(H))^2\V_H(\widehat{c}_n) 
\end{eqnarray}
for each $H\in\ ]0,1]$ (see Table \ref{tab2} below). In order to avoid confusion 
with the actual expectation and the actual variance, 
we call the expressions on the
right side of (\ref{asbias}) and (\ref{asvar}) the \textit{asymptotic expectation}
and the \textit{asymptotic variance}, respectively, of $\widehat{H}_n$.
See Figure \ref{fig1} (b), (c) for plots of 
the asymptotic bias, which is the difference between the asymptotic expectation and
$H$, and of the variance.
The computation of $\V_H(\widehat{c}_n)$
was carried out using formula (\ref{varapprox}) given below.

\begin{figure}[h]
\begin{picture}(0,100)
\put(0,0){\includegraphics[width=5cm]{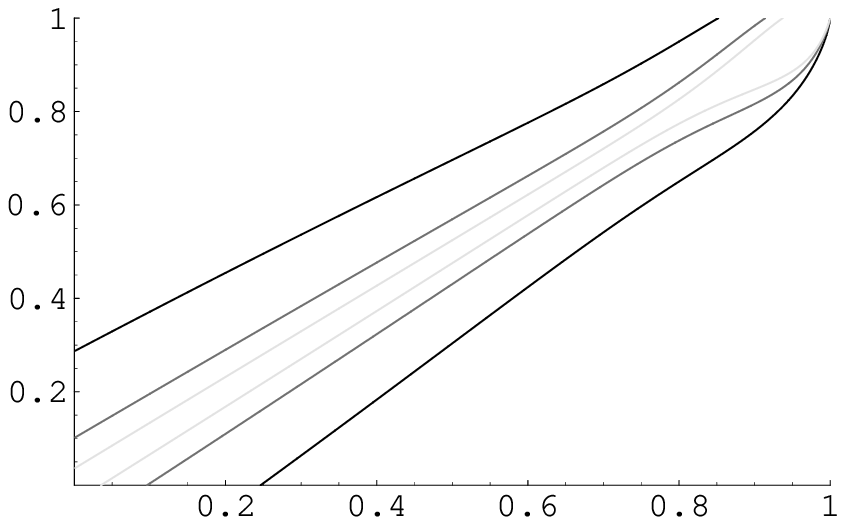}}
\put(70,90){(a)}
\put(160,0){\includegraphics[width=5cm]{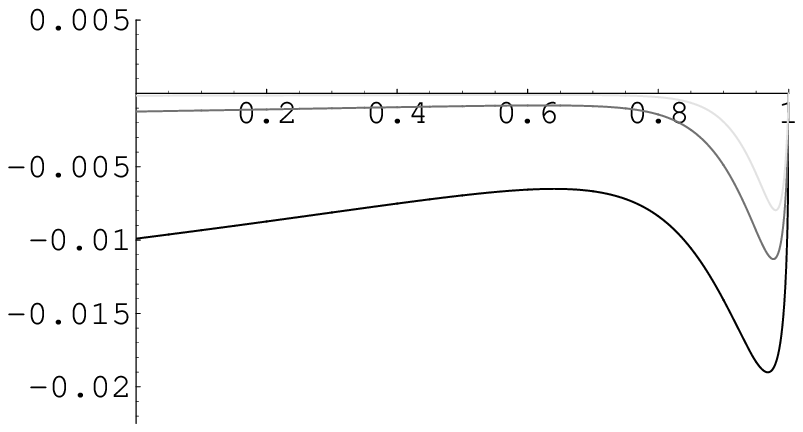}}
\put(230,90){(b)}
\put(320,0){\includegraphics[width=5cm]{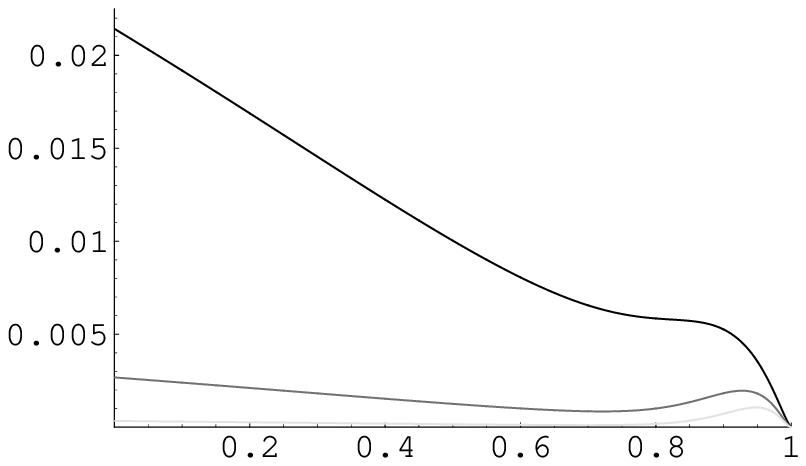}}
\put(390,90){(c)}
\end{picture}
\caption{{(a) 
Confidence intervals for the Hurst parameter obtained from $\widehat{H}_n$, for $n=128$ (black), $n=1024$ (dark grey)
and $n=8192$ (light grey). For the value $x$ on the abscissa, the ordinate shows the
confidence interval $K_n$
obtained if $\widehat{H}_n(\omega)=x$.
(b) Asymptotic bias of $\widehat{H}_n$, plotted against $H$. (c) Asymptotic variance of $\widehat{H}_n$, plotted against $H$.}}\label{fig1}
\end{figure}

\paragraph{Practical approximations.}
For the practical purpose of computing approximate numerical values
of $\V_H(\widehat{c}_n)$ for given $H\in\ ]0,1[$ and $n\in\N$,
first note that
the approximation of $\gamma_H(k)$ by the asymptotically equivalent expression
\begin{eqnarray*}
\widetilde{\gamma}_H(k) &:=& \frac{2(1-\rho_H(1))}{\pi^2(1+\rho_H(1))} H^2(2H-1)^2k^{4H-4}
\end{eqnarray*}
obtained in
Theorem \ref{ascov} can be improved in the following way.
For $(r_1,r_2)\in\R^2$
such that $(r_1,r_2,r_2,r_2)\in\cR$ (with $\cR$ as given in the preceding section),
let
\begin{eqnarray*}
F((r_1,r_2)) &:=& \Phi\circ\Sigma((r_1,r_2,r_2,r_2)).
\end{eqnarray*}
Further, define
\begin{eqnarray*}
\widetilde{\gamma}_H^{(m)}(k) &:=&
4\,\sum_{l=1}^m \frac{1}{(2l)!}
\frac{\partial^{2l} F}{\partial^{2l} r_2}\big((\rho_H(1),0)\big) \big(H(2H-1)k^{2H-2}\big)^{2l}
\end{eqnarray*}
for $H\in\ ]0,1[$, $m\in\N$ and $k\geq2$. Note that
$\widetilde{\gamma}_H^{(1)}(k)=\widetilde{\gamma}_H(k)$ for all $H\in\ ]0,1[$ and $k\geq2$.
For $m=3$, $\epsilon=0.01,0.001$ and
several values of $H$, Table \ref{tab1} provides the number
\begin{eqnarray*}
k_H^{(m)}(\epsilon) &:=& \min \big\{k\geq 2: \  \frac{|\widetilde{\gamma}_H^{(m)}(k)-\gamma_H(k)|}{\gamma_H(k)}<\epsilon  \big\},
\end{eqnarray*}
that is the least $k\geq 2$ such that the relative error caused
by approximating $\gamma_H(k)$ by $\widetilde{\gamma}_H^{(m)}(k)$ is less than $\epsilon$.
Because $\gamma_H(k)$ and $\widetilde{\gamma}_H^{(m)}(k)$ are asymptotically equivalent for all $m\in\N$,
the relative error tends to $0$ as $k$ tends to $\infty$ (see \cite{Lehmann99}), Lemma 1.1.1), that is,
$k_H^{(m)}(\epsilon)<\infty$ for every $\epsilon>0$.
Numerical experiments suggest that
$\widetilde{\gamma}_H^{(m)}(k)\leq\gamma_H(k)$ for all $m\in\N$,
$H\in\ ]0,1[$ and $k\geq2$.
As one can see in Table \ref{tab1}, $\widetilde{\gamma}_H^{(3)}(k)$ already
yields good approximations of $\gamma_H(k)$ for small $k$ if $H$ is not too close
to $1$. 

\begin{center}
\begin{table}[h]
\centering
\begin{tabular}{lllllllllll}
$H$ & $0.05$ & $0.15$ & $0.25$ & $0.35$ & $0.45$ & $0.55$ & $0.65$ & $0.75$ & $0.85$ & $0.95$ \\ \hline \vspace{-3mm} &&&&&&
\\\vspace{2mm}
$k_H^{(3)}(0.01)$ & $18$ & $16$ & $14$ & $12$  & $11$ & $9$ & $7$ & $5$ & $5$ & $226$\\
$k_H^{(3)}(0.001)$& $55$ & $50$ & $44$ & $38$  & $32$ & $26$ & $21$ & $15$ & $13$ & $10\,040$
\end{tabular}
\caption{The numbers $ k_H^{(3)}(\epsilon) $ for $\epsilon=0.01,0.001$ and several
values of the Hurst parameter.}
    \label{tab1}
\end{table}
\end{center}

By the formulas for the first partial
derivatives of $(\Phi\circ\Sigma)$ given in (\ref{firstpartder}) one obtains
\begin{eqnarray*}
\frac{\partial F}{\partial r_2}(\defn{s}) &=& 
\frac{2}{\pi(1-s_2^2)^{\frac{1}{2}}}\left(
\frac{1}{4}-\frac{1}{2\pi}\arcsin\frac{s_2(s_1-1)}{1+s_1-2s_2^2}
  \right)
\end{eqnarray*}
for every $\defn{s}=(s_1,s_2)$ in the domain of $F$. The
partial derivatives with respect to $r_2$ of
higher orders are obtained by easy calculations. Those of $4$th and $6$th order in $(\rho_H(1),0)$ 
are given by
\begin{eqnarray*}
\frac{\partial^4 F}{\partial^4 r_2}((\rho_H(1),0)) &=& \frac{4(1-\rho_H(1))(2+\rho_H(1))^2}{\pi^2(1+\rho_H(1))^3}, \\
\frac{\partial^6 F}{\partial^6 r_2}((\rho_H(1),0)) &=& \frac{16(1-\rho_H(1))(7+6\rho_H(1)+2\rho_H(1)^2)^2}{\pi^2(1+\rho_H(1))^5}.
\end{eqnarray*}
For a practical approximation of $V_H(\widehat{c}_n)$ we propose 
to sum up the precise numerical values $\gamma_H(k)$ obtained by (\ref{numform}) 
for $k=2,3,\ldots$ up to some specified 
$\widetilde{n}\leq n-1$,
and then to take the approximate values $\widetilde{\gamma}_H^{(3)}(k)$
for $k=\widetilde{n}+1, \ldots, n-1$, that is
\begin{eqnarray}\label{varapprox}
\V_H(\widehat{c}_n) \ \approx\
\frac{1}{n^2} \big( \, n\gamma_H(0) + 2\sum_{k=1}^{\widetilde{n}-1}(n - k)\,\gamma_H(k)
+ 2\sum_{k=\widetilde{n}}^{n-1}(n - k)\,\widetilde{\gamma}_H^{(3)}(k)\,
 \big).
\end{eqnarray}
For the calculations behind Figure \ref{fig1} we have chosen $\widetilde{n}$
depending on $H$, namely
$\widetilde{n}=\min\{k_H^{(3)}(0.01),250,n\}$. 

\section{Simulation studies}\label{sim}

For the simulation of fBm we used the algorithm
of Davies and Harte \cite{DaviesHarte87} and the random number generator of Mathematica 5.0.1.0.
For each of the values $H=0.55, 0.65, 0.75,$ $0.85, 0.95$ and the 
sample lengths $n=128, 1024, 8192$, we performed $50\,000$ simulations of fBm.
We only chose values $H>\frac{1}{2}$, because the increments of fBm
are particularly interesting
for the modelling of long range dependence in this case (see \cite{Beran94}).

\subsection{Performance of the ZC estimator}

Table \ref{tab2} shows simulations of $\E_H( \widehat{H}_n )$ and 
$\V_H( \widehat{H}_n )$, given
by the sample mean and the sample variance, respectively, of the $50\,000$ estimates
$\widehat{H}_n$ obtained for simulated samples of length $n$ with the true Hurst parameter $H$ . 
Additionally, Table \ref{tab2} provides the asymptotic expectation (abbreviated as.~exp.) and the
asymptotic variance (abbr. as.~var.) of $\widehat{H}_n$ as given by the expressions on the right side
of equations 
(\ref{asbias}) and (\ref{asvar}), respectively. As one can see,
the results obtained by the simulations and by equations 
(\ref{asbias}) and (\ref{asvar}), respectively,
agree very well. Only for $H=0.85, 0.95$ the
asymptotic expressions seem to yield
slightly larger values for $\V_H\big( \widehat{H}_n  \big)$
than the simulations.

Table \ref{tab2} also provides simulations of 
$\P_H( H \in K_n)$,
which is the probability that the true value of $H$ falls within
the interval $K_n$ obtained from formula (\ref{convintform}).
Theoretically, $\P_H( H \in K_n )$
tends to $0.95$ as $n$ tends to $\infty$ if $H<\frac{3}{4}$ 
(see (\ref{confidence1})). As one can see,
for $H=0.55, 0.65, 0.75$ the simulations of $\P_H( H \in K_n )$
always attain values larger than $0.95$. For $H=0.85$ this is still true
for $n=1024$ and $n=8096$, while the simulations of $\P_H( H \in K_n )$
are all clearly below $0.95$ for $H=0.95$. 

\begin{figure}[h]
\begin{picture}(0,100)
\put(0,0){\includegraphics[width=5cm]{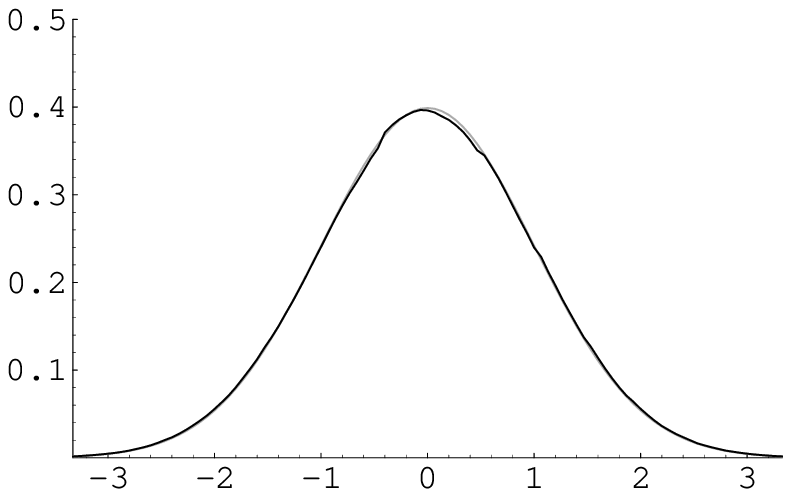}}
\put(70,90){(a)}
\put(160,0){\includegraphics[width=5cm]{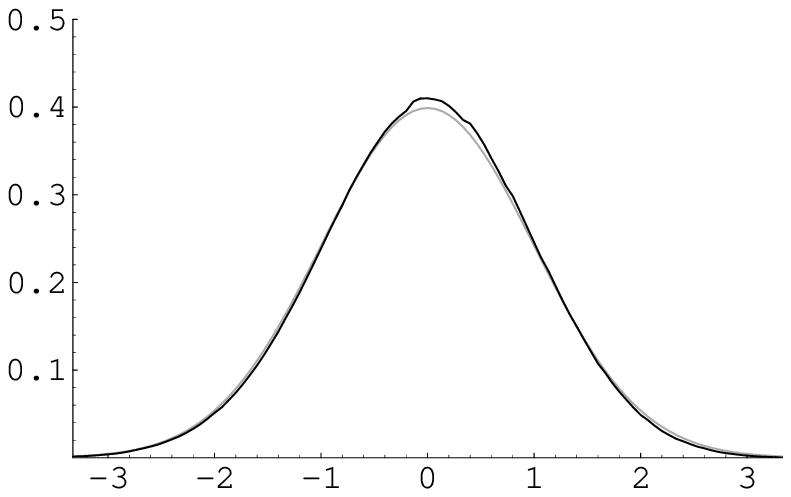}}
\put(230,90){(b)}
\put(320,0){\includegraphics[width=5cm]{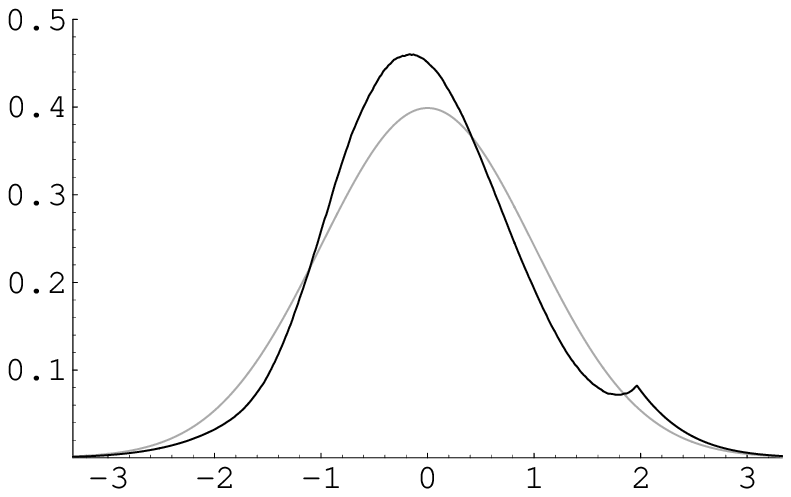}}
\put(390,90){(c)}
\end{picture}
\caption{{Standardized simulated distribution of $ \widehat{H}_n$ for $n=8192$
and (a) $H=0.55$ , (b) $H=0.75$  and (c) $H=0.95$, compared to the 
standard normal distribution (light grey).}}
\label{fig3}
\end{figure}

For $n=8192$ and $H=0.55,0.75,0.95$, Figure \ref{fig3} shows the standardized simulated distribution of $\widehat{H}_n$ with respect to $\P_H$
compared to the density of the standard normal distribution.
In accordance with (\ref{Hasnorm}),
for $H=0.55$ the standardized distribution of $\widehat{H}_n$ 
is very close to the standard normal distribution. For $H=0.75$ 
there is a slight difference between both curves, in particular,
the standardized distribution of $\widehat{H}_n$ seems to be more concentrated
around $0$ than the standard normal distribution. 
For $H=0.95$, Figure \ref{fig3} shows a remarkable contrast between the standardized distribution of $\widehat{H}_n$ and the standard normal distribution,
suggesting that $\widehat{H}_n $ is not asymptotically normally distributed
in this case.

\begin{center}
\begin{table}[ht]
\centering
\begin{tabular}{lllllll}
& $H$ &  $0.55$ & $0.65$ & $0.75$ & $0.85$ & $0.95$ \\\hline
\vspace{-3mm} &&&&&&
\\\vspace{2mm}
$n=128$ & $\E_H( \widehat{H}_n  )$ & $0.544$ & $0.643$ & $0.742$ & $0.839$  & $0.932$ \\\vspace{2mm}
& as. exp. & $0.543$ & $0.643$ & $0.743$ & $0.839$ & $0.932$ \\\vspace{2mm}
& $\V_H( \widehat{H}_n  )$ & $0.00935$ & $0.00754$ & $0.00616$ & $0.00543$  & $0.00389$ \\\vspace{2mm}
& as. var. & $0.009$ & $0.00722$ & $0.00609$ & $0.00572$ & $0.00354$ \\\vspace{2mm}
& $\P_H( H \in K_n )$ & $0.955$ & $0.96$ & $0.964$ & $0.931$ & $0.749$ 
\\\hline\vspace{-3mm} &&&&&&
\\\vspace{2mm}
$n=1024$ & $\E_H( \widehat{H}_n  )$ & $0.549$ & $0.649$ & $0.749$ & $0.848$  & $0.941$ \\\vspace{2mm}
& as. exp. & $0.549$ & $0.649$ & $0.749$ & $0.848$ & $0.941$ \\\vspace{2mm}
& $\V_H( \widehat{H}_n  )$ & $0.00113$ & $0.000912$ & $0.000849$ & $0.0012$ &
$0.00149$  \\\vspace{2mm}
& as. var. & $0.00113$ & $0.000913$ & $0.000863$ & $0.00134$ & $0.00181$ \\\vspace{2mm}
& $\P_H( H \in K_n )$ & $0.952$ & $0.952$ & $0.96$ & $0.954$ & $0.823$
\\\hline\vspace{-3mm} &&&&&&
\\\vspace{2mm}
$n=8192$ & $\E_H( \widehat{H}_n  )$ & $0.55$ & $0.65$ & $0.75$ & $0.849$ & $0.945$ \\\vspace{2mm}
& as. exp. & $0.55$ & $0.65$ & $0.75$ & $0.849$ & $0.945$ \\\vspace{2mm}
& $\V_H( \widehat{H}_n  )$ & $0.000141$ & $0.000116$ & $0.000121$ & $0.000316$ & 
$0.000796$ \\\vspace{2mm}
& as. var. & $0.000141$ & $0.000115$ & $0.000121$ & $0.000347$ & $0.00106$ \\\vspace{2mm}
& $\P_H( H \in K_n )$ & $0.95$ & $0.951$ & $0.953$ & $0.971$ & $0.884$ \\\hline
\end{tabular}
\caption{Simulations of $\E_H( \widehat{H}_n )$, 
$\V_H( \widehat{H}_n )$ and $\P_H( H \in K_n )$,
obtained each for $50\,000$ simulations of fBm. The rows `as. exp.' and `as. var.' show
the asymptotic expectation and variance, respectively, given by the right side of equations
(\ref{asbias}) and (\ref{asvar}), respectively. Note that,
according to (\ref{confidence1}), if $H<\frac{3}{4}$ then
$\P_H( H \in K_n)$ tends to $0.95$ as $n$ tends to $\infty$.}
    \label{tab2}
\end{table}
\end{center}

\subsection{Comparison to the HEAF estimator}
Let $(Y_k)_{k\in\N}$ be the increment process of an equidistant discretization of fBm,
and define $\overline{Y}_n:=\frac{1}{n}\sum_{k=1}^n Y_k$ for $n\in\N$.
The HEAF estimator (HEAF meaning `Hurst Exponent from Autocorrelation Function') proposed
by Kettani and Gubner \cite{KettaniGubner02} puts the estimate
\begin{eqnarray*}
\widehat{\rho}_n &=& \frac{\sum_{k=1}^{n-1}(Y_k-\overline{Y}_n)(Y_{k+1}-\overline{Y}_n)}
{\sum_{k=1}^n(Y_k-\overline{Y}_n)^2},
\end{eqnarray*}
of $\rho_H(1)$ into the monotonic functional
relation $\rho_H(1) = 2^{2H-1}-1$ obtained from (\ref{hdef}). In order to receive only
finite non-negative estimates of $H$, let
\begin{eqnarray*}
\widehat{H}_n^{\text{HEAF}} &=& \frac{1}{2}\big( 1 + 
\log_2\big(1+  \max\{ -1/2,\widehat{\rho}_n \}   \big) \big)
\end{eqnarray*}
for $n\in\N$. Table \ref{tab3} shows simulations of $\E_H( \widehat{H}_n^{\text{HEAF}})$
and $\V_H( \widehat{H}_n^{\text{HEAF}} )$, each
obtained for $50\,000$ simulations of fBm.

As one can see, $\widehat{H}_n^{\text{HEAF}}$ has smaller variance
but larger bias than the ZC estimator. For $n=128$ the variance
of the ZC estimator is about $2$ times larger than the variance of $\widehat{H}_n^{\text{HEAF}}$,
for $n=1024$ it is about $2$ to $4$ times larger, and for $n=8192$
it is about $2$ to $8$ times larger. Notice that the highest ratios are obtained
for $H=0.95$.

The loss of efficiency of $\widehat{H}_n$
relative to $\widehat{H}_n^{\text{HEAF}}$ is not surprising. 
In particular, $\widehat{H}_n$
can be seen as the estimator obtained by plugging the estimate 
$\sin(\pi(1/2-\widehat{c}_n))$
of $\rho_H(1)$ into the monotonic functional relation
$\rho_H(1) = 2^{2H-1}-1$ (compare to (\ref{proofd=2}) and
(\ref{HfromZeta}), respectively),
and clearly $\sin(\pi(1/2-\widehat{c}_n))$ loses efficiency relative
to the estimate $\widehat{\rho}_n$ of $\rho_H(1)$.
However, because of the slow convergence of 
the sample mean for increments of fBm with Hurst parameter $H>\frac{1}{2}$ (see \cite{Beran94}), $\widehat{H}_n^{\text{HEAF}}$
has a much larger bias than the ZC estimator, in particular for large values of $H$.

In a concluding remark, the ZC estimator of the Hurst parameter performs comparatively well,
although it does not regard the metric structure
but only the order relations between values of realizations.
Notice that generally, for fixed $r\in S_d$ with $d\in\N$,
the computation of $\widehat{p_{\overline{r},}}_n$
can be carried out by 
asymptotically 
$4\,(\sharp\,\overline{r})\,d\,n$ computational steps (see \cite{KellerLaufferSinn05}).

\begin{center}
\begin{table}[ht]
\centering
\begin{tabular}{lllllll}
& $H$ &  $0.55$ & $0.65$ & $0.75$ & $0.85$ & $0.95$ \\\hline
\vspace{-3mm} &&&&&&
\\\vspace{2mm}
$n=128$ & $\E_H( \widehat{H}_n^{\text{HEAF}}  )$ & $0.538$ & $0.628$ & $0.712$ & $0.787$  & $0.849$ \\\vspace{2mm}
& $\V_H( \widehat{H}_n^{\text{HEAF}}  )$ & $0.00377$ & $0.00322$ & $0.00267$ & $0.00218$  & $0.00167$
\\\hline\vspace{-3mm} &&&&&&
\\\vspace{2mm}
$n=1024$ & $\E_H( \widehat{H}_n^{\text{HEAF}}  )$ & $0.548$ & $0.646$ & $0.739$ & $0.824$  & $0.893$ \\\vspace{2mm}
& $\V_H( \widehat{H}_n^{\text{HEAF}}  )$ & $0.000468$ & $0.000399$ & $0.000378$ & $0.000373$ & $0.000328$ 
\\\hline\vspace{-3mm} &&&&&&
\\\vspace{2mm}
$n=8192$ & $\E_H( \widehat{H}_n^{\text{HEAF}}  )$ & $0.55$ & $0.649$ & $0.746$ & $0.837$ & 
$0.912$ \\\vspace{2mm}
& $\V_H( \widehat{H}_n^{\text{HEAF}}  )$ & $0.0000584$ & $0.0000515$ & $0.0000573$ & $0.0000836$ &
$0.000101$ \\\hline
\end{tabular}
\caption{Simulations of the mean and of the variance of the HEAF estimator.}
    \label{tab3}
\end{table}
\end{center}

\section*{Acknowledgement}
We would like to thank Lutz Mattner for many fruitful discussions and
important hints concerning Theorem \ref{RaoBlackwell} and Theorem \ref{ascov}.


\begin{thebibliography}{1}

\bibitem{Arcones94} Arcones,~M.~A., Limit theorems for nonlinear functionals of a stationary Gaussian sequence of vectors. \textit{Annals of Probability} 22 (1994), 2242-74.

\bibitem{Bandt04} Bandt,~C., Ordinal time series analysis. \textit{Ecological Modelling} 182 (2005),
229-38.

\bibitem{BandtShiha05} Bandt,~C. and Shiha,~F.,
Order patterns in time series. \textit{J. Time Ser. Anal.} 28 (2007), 646-65.

\bibitem{Beran94} Beran,~J., \textit{Statistics for Long-Memory Processes}. 
London: Chapman and Hall (1994).

\bibitem{Cheng69} Cheng,~M.~C., The orthant probability of four Gaussian variates. \textit{Ann. Math. Statist.} 40 (1969), 152-61.

\bibitem{Coeurjolly00}
Coeurjolly,~J.~F., Simulation and identification of the fractional
Brownian motion: A bibliographical and comparative study. \textit{J.
Stat. Software} 5 (2000).

\bibitem{Cornfeld82}
Cornfeld,~I.~P., Fomin,~S.~V.~ and Sinai,~Ya.~G., \textit{Ergodic Theory}. 
Berlin: Springer-Verlag (1982).

\bibitem{DaviesHarte87} Davies,~R.~B. and Harte,~D.~S., 
Tests for Hurst effect. \textit{Biometrika} 74 (1987), 95-101.

\bibitem{EmbrechtsMaejima02} Embrechts,~P. and Maejima,~M., \textit{Selfsimilar Processes}. 
New Jersey: Princeton University Press (2002).

\bibitem{HoSun87} Ho,~H.-C. and Sun,~T.~C., A central limit theorem for noninstantaneous filters of a stationary Gaussian process.
\textit{J.~Multivariate Anal.} 22 (1987), 144-55.

\bibitem{Kedem94} Kedem,~B., \textit{Time Series Analysis by Higher Order Crossings}.
New York: IEEE Press (1994).

\bibitem{KellerLaufferSinn05} Keller,~K., Lauffer,~H. and Sinn,~M.,
Ordinal analysis of EEG time series.
\textit{Chaos and Complexity Letters} 2 (2005), 247-58.

\bibitem{KellerSinnEmonds07} Keller,~K., Sinn,~M. and Emonds,~J., 
Time series from the ordinal viewpoint. \textit{Stochastics and Dynamics} 2 (2007), 247-72.

\bibitem{KettaniGubner02} Kettani,~H. and Gubner,~J.~A., 
A novel approach to the estimation of the Hurst parameter in self-similar traffic.
\textit{IEEE Trans. Circuits Syst. II} 53 (2006), 463-67.

\bibitem{Lehmann86} Lehmann,~E., 
\textit{Testing Statistical Hypotheses, Second Edition}. New York: Springer-Verlag (1986).

\bibitem{Lehmann99} Lehmann,~E., 
\textit{Elements of Large Sample Theory}. New York: Springer-Verlag (1999).

\bibitem{MarkovicKoch05} Markovi\'c,~D. and Koch,~M., 
Sensitivity of Hurst parameter estimation to periodic signals in time series and filtering approaches.
\textit{Geophysical Research Letters} 32 (2005), L17401.

\bibitem{Pfanzagl94}
Pfanzagl,~J., \textit{Parametric Statistical Theory}. Berlin, New York: de Gruyter (1994).

\bibitem{Plackett54} Plackett,~R.~L., 
A reduction formula for normal multivariate integrals. \textit{Biometrika} 41 (1954), 351-60.

\bibitem{Shietal05} Shi,~B., Vidakovic,~B., Katul,~G. and Albertson,~J.~D.,
Assessing the effects of atmospheric stability on the fine structure of
surface layer turbulence using local and global multiscale
approaches. \textit{Physics of Fluids} 17 (2005), 055104.

\end{thebibliography}
\end{document}